\documentclass[12pt, leqno]{amsart}

\usepackage{amsfonts,amsmath,amsthm,amssymb,bbm}
\usepackage{latexsym}
\usepackage{enumerate}
\usepackage[usenames]{color}
\usepackage{hyperref,latexsym,multicol}
\usepackage[usestackEOL]{stackengine}
\usepackage[pagewise]{lineno}
\newcommand{\defc}{black}

\usepackage{xcolor}
\newtheorem{theorem}{Theorem}[section]
\newtheorem{conj}{Conjecture}
\newtheorem*{theoremBS}{Theorem A.}
\newtheorem{lemma}[theorem]{Lemma}

\theoremstyle{definition}

\usepackage[ansinew]{inputenc}
\usepackage{ifthen}
\usepackage{animate}

\newcommand{\al}{\alpha}

\newcommand{\ga}{\gamma}
\newcommand{\la}{\lambda}
\newcommand{\La}{\Lambda}

\newcommand{\ve}{\varepsilon}
\newcommand{\Z}{\mbox{$\mathbb Z$}}

\newcommand{\be}{\begin{equation}}
\newcommand{\ee}{\end{equation}}

\title[Erd\H{os}-Selfridge superelliptic curve]{Rational solutions to the Variants of Erd\H{o}s- Selfridge superelliptic curves}

\author[Das]{Pranabesh Das}
\address{Department of Mathematics, Xavier University of Louisiana,  1 Drexel Dr, New Orleans, LA 70125, USA}
\email {pranabesh.math@gmail.com}

\author[Laishram]{Shanta Laishram}
\address{Stat Math Unit, 
Indian Statistical Institute, 7 SJSS Marg, New Delhi, 110016, India}
\email{shanta@isid.ac.in}

\author[Saradha]{N. Saradha}
\address{INSA Senior Scientist, 
DAE-Centre for Excellence in Basic Sciences, 
University of Mumbai, Mumbai, 400098, India}
\email{saradha54@gmail.com}

\author[Sharma]{Divyum Sharma}
\address{Department of Mathematics, Birla Institute of Technology and Science, Pilani, Rajasthan,  333031, India}
\email{divyum.sharma@pilani.bits-pilani.ac.in}

\date{\today}

\begin{document}

\begin{abstract}
For the superelliptic curves of the form
$$ (x+1) \cdots(x+i-1)(x+i+1)\cdots (x+k)=y^\ell$$
with $x,y \in \mathbb{Q}$, $y\neq 0$, $k \geq 3$, $1\leq i\leq k$, $\ell \geq 2,$ a prime,
Das, Laishram, Saradha and Edis showed that the superelliptic curve has no rational points for $\ell\geq e^{3^k}$. In fact the double exponential bound, obtained in these papers is far from the reality. In this paper, we study the superelliptic curves for small values of $k$. In particular, we explicitly solve the above equation for $4 \leq k \leq 8.$
\end{abstract}

\maketitle

\section{\bf Introduction}
In a remarkable work in 1975,  Erd\H{o}s and Selfridge\cite{ES} showed that the Diophantine equation 
\begin{equation}\label{er-sel}
(x+1)(x+2)\cdots(x+k)= y^{\ell}
\end{equation}
has no solutions in positive integers $x,y,k,$ and $\ell$ with $k,\ell \geq 2.$\\

Since their result, several
variations of the equation have been considered by many mathematicians and integral solutions have been investigated. The equation \eqref{er-sel} has been studied in a more general setting of arithmetic progressions. We refer to the survey articles of Shorey \cite{Sh1} and \cite{Sh2} for various results in this direction.\\

We state the following conjecture on the generalizations of equation \eqref{er-sel} widely attributed to Erd\H{o}s.

\begin{conj}\label{conj1}(Erd\H{o}s)
There is a constant $k_0$ such that the Diophantine equation
\begin{equation}\label{er-conj}
(x+d)(x+2d)\cdots(x+kd)= y^{\ell}, \ \gcd(x,d)=1
\end{equation}
has no solutions in positive integers $x, y, k, d, y, \ell$ with $\ell \geq 2$ and $k\geq k_0$. That is,
if the equation has a solution in $x, |y|>1$, then $k$ is bounded by an absolute constant.
\end{conj}
The natural condition $\gcd(x,d)=1$ is imposed to avoid trivial solutions. We note the equation \eqref{er-conj} has infinitely many non-trivial solutions for $(k,\ell)=(3,2).$
The above conjecture is very difficult and still very much out of the reach. Currently, the Erd\H{o}s conjecture has been verified unconditionally only subject to a variety of additional assumptions such as $d$ is fixed (Marzalek \cite{Mar}) or $\ell$ and $\omega(d)$(the number of distinct prime divisors of $d$) both are fixed (Shorey and Tijdeman \cite{ShTi}) or $P(d)$ \textcolor{\defc}{(the greatest prime factor of $d$)} is fixed and $\ell \geq 3$ (Shorey \cite{Sh1988}) or $x$ is fixed and $\ell \geq 7$ (Shorey \cite{Sh1996}).

Recently, in an elegant paper, Bennett and Siksek \cite{BeSi2020} proved a weakened version of Conjecture \ref{conj1}. They proved the following.
\begin{theoremBS}[Bennett-Siksek\cite{BeSi2020}]
There is an effectively computable absolute constant $k_0$ such that if $k\geq k_0$ is a positive integer, then any solution in integers to equation \eqref{er-conj} with prime
exponent $\ell$ satisfies either $y = 0$ or $d = 0$ or $\ell \leq exp(10k)$.
\end{theoremBS}

In a different direction, rational points on the equation \eqref{er-conj} has been studied in the literature. Next, we introduce some notation which will be used throughout the paper.\\

Let $x\in \mathbb Q$ and $k \geq 2$ be an integer. For any integer $n\geq 1,$
let $P(n)$ denote the greatest prime factor of $n$ and take $P(1)=1.$ 
Further, let $\nu_p(n)$ denote the power of prime $p$ in the prime factorization of $n$ with $\nu_p(1)=0.$ 
Put
$$ \Delta_0=(x+1)\cdots (x+k)$$
and for $1 \leq i \leq k,$ let
$$\Delta_i=(x+1)\cdots (x+i-1)(x+i+1)\cdots (x+k).$$
Bennett and Siksek \cite{BeSi} considered rational solutions of
\be\label{ES-curve}
\Delta_0=y^\ell
\ee
in $x$ and $y$ with \textcolor{\defc}{$y\neq 0$} and $\ell \geq 2,$ a prime. They showed that if \eqref{ES-curve} holds, then
\be\label{upperbound}
\ell \leq e^{3^k}.
\ee
We note that the equation\eqref{ES-curve} defines a superelliptic curve of genus at least $(\ell-1)(k-1)/2.$ Since for a fixed pair $(k,\ell)$ except $(k,\ell)=\{(2,2),(2,3),(3,2)\}$ the superelliptic curve has genus $>1$, hence by a theorem of Falting it follows from Bennett and Siksek \cite{BeSi} that, for $\ell> 4$, the number of rational points on the superelliptic curve \eqref{ES-curve} is finite. Although, conjecturally the curve \eqref{ES-curve} has no non-trivial rational points for $\ell \geq 3$. The conjecture is due to Sander \cite{San} and it states that 
\begin{conj}[Sander\cite{San}]\label{conj2}
The superelliptic curve \eqref{ES-curve} has no rational points for $\ell \geq 3$ with $y\neq 0.$
\end{conj}
We are still quite far away from the conjecture. 
There has been some progress on Conjecture \ref{conj2}, for small values of $k$. Sander \cite{San} and Lakhal and Sander \cite{LaSan} studied \eqref{ES-curve} for $k\leq 6.$ In \cite{BBGH}, Bennett, Bruin, Gy{\H ory}, and Hajdu, solved \eqref{ES-curve} completely for $k \leq 11.$ This was further extended to $k\leq 34$ in \cite{GHP} by Gy{\H ory}, Hajdu and Pint\'{e}r.\\ 

The above problems can also be considered as a rational analogue to the Schinzel-Tijdeman theorem on integral solutions to the superelliptic equation
$
f(x)=y^l
$
where $f(x)$ is a polynomial.\\

In recent papers, Das, Laishram, Saradha \cite{DLS} and Edis\cite{S} extended the result of \cite{BeSi}  to the equation
\begin{equation}\label{ES-oneterm}
\Delta_i=y^\ell, 1 \leq i \leq k, y\neq 0.
\end{equation}
by showing that \eqref{upperbound} is valid. See also Saradha\cite{Sar1} when more than one term is omitted in $\Delta_0.$ This double exponential bound is very large and it is desirable to improve the bound. Towards this, Saradha \cite{Sar2} showed that if \eqref{ES-curve} has a positive rational solution, then either the denominator of the solution is large or $\ell\leq k.$\\

We notice that \eqref{ES-oneterm} has solutions for $x=-j; y=0$ for $1\leq j \leq k$ and $i\neq j$. We call these solutions as \textit{trivial solutions}.\\

\textcolor{\defc}{Our goal in this paper to study the equation \eqref{ES-oneterm} for small values of $k.$ It is difficult to solve equation \eqref{ES-oneterm} explicitly even for small values $k$. In this paper, we explicitly find the non-trivial rational points of the curve \eqref{ES-oneterm} for $4\leq k\leq 8.$ The ideas used in this paper, in principle may work for larger values of $k$ but the subsequent combinatorial arguments would be very involved for $k\geq 10$.} \\

When $k=3$, Shen and Cai \cite{ShCa} showed the following result. See \cite{ShCa}[Corollaries 1.3 and 1.4].

{\it Let $k=3.$ Then the only non-trivial rational points on \eqref{ES-oneterm} with $i=2$ are given by
$$(x,y)=
\begin{cases}
\{(-2,-1),(-5,2),(1,2)\}\ {\rm if}\ \ell=3,\\
(-2,-1)\ {\rm if}\ \ell\geq 5.
\end{cases}
$$
For $\ell=2,$ we have
$$(x,y)=\bigg(\frac{3c_1^2-c_2^2}{c_2^2-c_1^2},\frac{2c_1c_2}{c_2^2-c_1^2}\bigg)$$
with co-prime integers $c_1\neq \pm c_2.$}
\vskip 2mm
  
In this paper, we consider \eqref{ES-oneterm} for $ 4\leq k\leq 8.$

\begin{theorem}\label{Th1-3-8}
Assume \eqref{ES-oneterm} for $4\leq k \leq 8$ and $(k,
\ell,i)\neq (4,3,2), (7,3,4).$ Then there are no solutions except when
\begin{enumerate}[(i)]
\item  $(k,\ell,i)= (7,2,2), (7,2,6)$, in which case the only non-trivial rational solutions are given by\\
$$(x,y)=\left(\frac{-37}{7}, \pm \frac{720}{7^3}\right), \left(\frac{-19}{7}, \pm \frac{720}{7^3}\right),$$ respectively.
\item $(k,\ell,i)= (5,2,2), (5,2,4)$ in which case
there exist infinitely many solutions corresponding to the non-trivial rational points of the elliptic curve $F: y^2= x^3 + 8x^2 + 12x.$
\end{enumerate}
In particular, there are no solutions for $\ell>3.$\\
\end{theorem}
\noindent When $(k,\ell,i)=(5,2,2),$ some of the non-trivial solutions are given by $(x,y)\in $
$$\left\{ (-7,\pm 12), \Big(\frac{-11}{3},\frac{\pm 8}{9}\Big), \Big(\frac{-17}{5},\frac{\pm 24}{25} \Big), \Big(\frac{-3}{7}, \frac{\pm 240}{49}\Big), \Big(\frac{-119}{23},\frac{720}{23^2}\Big)\right\}.$$\\
When $(k,\ell,i)=(7,3,4)$, our elimination techniques fail. This is due to the fact that all the corresponding elliptic curves we used for elimination are of rank $1$ or more. In this case, there exists a non-trivial solution given by
$$(x,y)=\left(\frac{-17}{7},  \frac{120}{7^2}\right)$$
which lead to an obstruction for our elimination techniques. Although, we expect that the given solution is the only rational solution. Additionally, we could not apply the Chabauty arguments of Bruin and Stoll directly since it is not a hyperelliptic curve for the case $(k,\ell,i)=(7,3,4)$. It would be interesting to study the rational points separately for the superelliptic curve for $$(x+1)(x+2)(x+3)(x+5)(x+6)(x+7)=y^3$$.\\

\noindent When $(k, \ell)=(4,3),$ the solutions to \eqref{ES-oneterm} will arise from the integral solution to the equation $u^3+ 2v^3=3w^3$ in coprime integers $u,v,w.$ These solution can be obtained from the rational points of the curve $U^3 + 2V^3= 3$. The curve $U^3 + 2V^3= 3$ has infinitely rational solutions appearing from the bi-rationally equivalent Weierstrass curve $Y^2= X^3 - 432\times 6^2$ of rank $1$ under the mutual inverse transformations $U=-\frac{24\theta^2X}{Y-216}, V=\frac{Y+216}{Y-216}\theta$ where $\theta=\sqrt[3]{\frac{3}{2}}$.  This may lead to infinitely many solutions in co-prime integers $u,v,w$.\\

As a prelude to the proof of Theorem \ref{Th1-3-8}, we show the following result.  

\begin{theorem}\label{Th2}
Suppose \eqref{ES-oneterm} has a rational solution for $4\leq k \leq 8$. 
Then $\gcd(k-1,\ell)>1.$ 
\end{theorem}
In Sections 2 and 3, we give required preliminaries and several lemmas on Fermat    equations. 
In Section 4, we consider the cases of $\ell>2$ with $\gcd(k-1, \ell)=1.$ 
 In Sections \ref{gcd>1}-$9$, we consider the case $\gcd(k-1, \ell)=\ell$ and finally in Section \ref{PfTh1}, 
 we consider $\ell=2$ and complete the proofs of Theorems \ref{Th2} and \ref{Th1-3-8}. 
\vskip 5mm

A main task is to form suitable ternary equations of signature $(\ell,\ell,\ell)$ or $(\ell,\ell,2),\ell>3$, and apply modularity methods to show that these equations have no non-trivial solution. Although this method is already used in \cite{GHP}, in our case the number of equations to be considered is {\it many},  since each equation  in \eqref{ES-oneterm} gives rise to $\bigg\lceil \frac{k}{2} \bigg\rceil-1$ equations as $i$ varies.  Further, one needs to use combinatorial arguments and the distribution of the small primes $2,3,5,7$ among the terms in the product $\Delta_i.$ Several new ternary equations of the signature $(7,7,2)$ are shown to have no non-trivial integral solutions in Lemma 2.9. For $\ell=3,$ we use an old result of Selmer \cite{Sel51}. The case $\ell=2$ requires several elliptic curves to be shown to have rank 0, and many quartic curves to have no non-trivial solution. Further the signs of the terms in the product $\Delta_i$ also play a crucial role.

\section{\bf Lemmas on ternary equations}

We begin this section with an old result of Selmer \cite[Tables $2^a,4^a,4^b$ and $6$]{Sel51}
on cubic ternary equations $Ax^3+By^3+Cz^3=0.$ We look for non-trivial solutions i.e. with $xyz\neq 0.$

\begin{lemma}\label{ell3}
Let $a$ and $ b$ be pairwise co-prime positive integers. Then the equation
$$x^3+ay^3+bz^3=0$$
has no solution in non-zero integers $x,y,z$ with $\gcd(x,y,z)=1$ for  
\begin{equation*}
ab \in \{3, 4, 5, 10,18,25, 36, 45, 60,100,150, 225,300\}. 
\end{equation*}
For $ab=2$, the equation has the solution $(1,1,-1)$ or $(1,-1,1)$.
\end{lemma}

For $\ell=5,$ we use the following lemma due to Kraus (see \cite[Proposition 6.1]{BBGH}). 
\begin{lemma}\label{ell5}
Let $c$ be a positive integer with $P(c)\leq 5.$ If the Diophantine equation
$$x^5+y^5=cz^5$$
has solutions in non-zero co-prime integers $x,y$ and $z,$ then $c=2$ and $x=y=z=\pm  1.$
\end{lemma}

\begin{lemma}\label{eq1}
 Let $\ell \geq 3$ and $\alpha \geq 0$ be integers. Then the equation
 \begin{align*}
 x^\ell+y^\ell=p^\alpha z^\ell \quad {\rm where} \quad \ p\in\{2,3\}\ and\ (\ell,p,\alpha)\neq (3,3,2)
 \end{align*}
 in relatively prime integers $x, y, z$ with $xyz\neq 0, 1$, has no solution. 
\end{lemma}
The results for $p=2$ were established by Wiles\cite{Wi} for
$\alpha=0,$ by Darmon and Merel \cite{DaMe}
for $\alpha=1,$ and by Ribet \cite{Ri} for
$\alpha> 1.$ Let $p=3.$ Suppose $\ell=3.$ The case $\alpha=0$ is classical.  The case $\alpha=1$
follows from Lemma \ref{ell3}. The result for all other cases were proved by Serre \cite{Serre}.
The next lemma is \cite[Lemma 13]{SaSh1}.
\begin{lemma}\label{eq2}
Let $\ell \ge 5$.  Let $a, b, c$ be pairwise co-prime positive integers with $abc\in \{2^u3^v, 2^u5^v\}$
where $u$ and $v$ are non-negative integers with $u\geq 4$. Then the equation 
\begin{align*}
ax^{\ell}+by^{\ell} = cz^{\ell}
 \end{align*}
 has no solution in pairwise co-prime non-zero integers $x, y$ and $z$.
\end{lemma}

The following lemma is due to Bennett \cite{Be} when $k=5, \ell\geq 7, P(b)=3$ and the remaining cases are covered in \cite{BBGH}.

\begin{lemma}\label{k5}
Let $m$ and $s$ be non-zero co-prime integers and $k\in \{3, 4, 5\}$. Then the equation 
\begin{align*}
(m+s)(m+2s)\cdots (m+ks)= by^{\ell}
\end{align*}
has no solutions in non-zero integers $b, y, \ell$ with $\ell\geq 2$ and $P(b)\leq P_{k, \ell}$ where 
$P_{k, \ell}$ is given by 
$$
\begin{matrix}
k & \ell=2 & \ell=3 & \ell=5 & \ell\geq 7\\
3 & - & 2 & 2& 2\\
4 & 2 & 3 & 2& 2\\
5 & 5 & 3 & 3& 3.
\end{matrix}
$$
\end{lemma}

The next result is a simple consequence of \cite{Be}, see also \cite[Proposition 2.3]{GHP}.

\begin{lemma}\label{GHP}
Let $\ell\ge 7$ be a prime. Then the equation 
\begin{align*}
(m + s)(m+2s)(m + 4s)(m + 5s) = by^\ell
\end{align*}
has only the solutions $(m, s, b, y)=(\pm  3, \mp 1, 4, 1)$ in non-zero 
integers $m, s, b, y$ with $\gcd(m, s) = 1$ and $P(b) \leq 3$.
\end{lemma}

The next lemma is part of \cite[Proposition 3.1]{BBGH}.

\begin{lemma}\label{eq4}
Let $\ell \ge 7$ be prime, $\al, \beta$ be non-negative integers and let $a, b$ be co-prime positive integers.  
Then the following equations have no solution in non-zero co-prime integers $(x,y,z)$ with $xy\neq  \pm 1$:
\begin{enumerate}
\item[(i)] $x^{\ell} + 2^\al y^{\ell} = 3^\beta z^2$, $\al\neq 1$. 
\item[(ii)] $x^{\ell} + 2^{\alpha} y^{\ell} = 3z^2$ with $p\mid xy$ for $p\in\{5,7\}$. 
\item[(iii)] $x^{\ell} + y^{\ell} = cz^2, $ $c\in \{1, 2, 3, 5, 6, 10\}$.  
\item[(iv)] $ax^{\ell} + by^{\ell} = z^2$, $P(ab) \le 3$ with $p\mid xy$ for $p\in\{5,7\}$. 
\item[(v)] $ax^{\ell} + by^{\ell} = z^2$, $P(ab) \le 5$ with $7|xy$ and $\ell\geq 11$. 
\end{enumerate}
\end{lemma}

For $\ell=7,$ we need to solve more cases which are not covered by the above lemmas.
The equations considered in Lemma \ref{ell7} are of the form $(\ell,\ell,2)$. A recipe for attaching the required Frey curves and their respective conductors, discriminant has been described in a paper by Bennett and Skinner \cite{BS}. The application of modular method for generalized Fermat equation of    $(\ell,\ell,2)$ has also been discussed. We follow the same strategy for equations in Lemma \ref{ell7}.
\begin{lemma}\label{ell7}
Let $\alpha, \beta, a, b$ be non-negative integers \textcolor{\defc}{with $\gcd(a,b)=1$}. Then the following equations have no solution in non-zero co-prime integers $(x,y,z)$ with $xy\neq  \pm  1$:
\begin{enumerate}[(i)]
\item $ax^7 + by^7 = z^2, ab\in \{2^{\alpha}3^{\beta}: \alpha \geq 6 \}$.
\item $ax^7+by^7= z^2,\\
 ab\in \{ 2^{\alpha}3^{\beta}7^{\delta}: (\alpha,\beta,\delta)\neq (1,\beta,\delta) \ {\rm for} \ \beta\delta \geq 1\}, 
5\mid {xy}.$
\item $x^7 + 3^{\beta}y^7 = 2z^2,  \ 5 \mid xy.$
\end{enumerate}
has no non-zero integer solution.
\end{lemma}
\begin{proof}
Let $f= \sum_{m=1}^{\infty} c_m q^m $ where $q=e^{2\pi iz}$ be a newform of weight $2$ with trivial Nebentypus character. We assume 
that the equation has a non-zero integer solution, we associate a Frey curve $E/\mathbb{Q}$ with 
corresponding mod $\ell$ Galois representation
$$\rho_{\ell}^{E}: Gal (\mathbb{Q}/\bar{\mathbb{Q}})\xrightarrow{} GL_2(\mathbb{F}_{\ell})$$
on the $\ell$-torsion $E[\ell]$ on $E.$ This representation arises from a cuspidal newform $f$ of weight $2$ with trivial Nebentypus character of level $N_n(E)$  by \cite[Lemma 3.2, 3.3]{BS}. We have \textcolor{\defc}{$n=7$ in our case, so we simply denote $N_7(E)$ by $N$.}\\

The respective level of the newforms that are associated to equation $(i)$ are $N\in \{ 1, 2, 3, 6 \}$ according to \cite[Lemma 2.1, 3.2]{BS}. This is a contradiction.

For equation $(ii)$, we need only to consider the corresponding Frey curve that will arise from the weight $2$ newforms with trivial characters of level 
\begin{align*}
N\in &\{14,21,24,32,42,56,64,84,96,128,168,192,224,384,
448,672,\\ &896,1344\}.
\end{align*}
Since  $5|xy$, the Frey curve has a 
multiplicative reduction at $5$. Therefore $7|\mathcal{N}_{K/Q}(c_5\pm  6)$ from \cite[Proposition 4.2]{BS} as $5$ is co-prime to $7N$. \textcolor{\defc}{We find this to be untrue by
checking for 
all the newforms of the corresponding levels.} Hence there is no solution. 

For equation $(iii)$, we find that $N\in\{256, 768\}.$ Similar to equation $(ii)$, we conclude that $7|\mathcal{N}_{K/Q}(c_5\pm  6)$. We have $$c_5=0, \pm 2,\pm   4, \pm  4\sqrt{2}, \pm  4\sqrt{3},$$ for weight $2$ newforms of level $256,768$ with trivial Nebentypus character, which is a contradiction.
\end{proof}

In the next lemma, we follow the paper by Halberstadt and Kraus\cite{HalKra} where strategies to tackle generalized Fermat equation of signature $(n,n,n)$ \textcolor{\defc}{have} been discussed. We refer to \cite[Section 15.8.1]{HC} for a detailed computation.
\begin{lemma}\label{n-n-n-lemma}
Let $A,B,C,X,Y,Z$ be non-zero integers such that $AX,\\BY, CZ$ are pairwise co-prime. Suppose that $29,43\mid XYZ$. Then the equation
\begin{equation*}
    AX^7 + BY^7 = CZ^7
\end{equation*}
has no solutions for $ABC=2^{\alpha}3^{\beta}5^{\ga}7^{\delta}$ for $\alpha \geq 4, \beta\gamma\delta >0 $ and $7\nmid \beta\gamma\delta.$
\end{lemma}

\begin{proof}
We begin with the assumption that $(x,y,z)$ is a solution to the equation. We use the recipes given by Halberstadt and Kraus\cite{HalKra} to find the suitable Frey curves for the Fermat type equations of signature $(\ell,\ell,\ell)$. We attach the Frey curve $E: Y^2= X(X-Ax^7)(X+By^7)$. We write $R=ABC.$ From our assumption we have $\text{Rad}_2(R)=105.$ We apply Ribet's level lowering. The Serre conductor of the Frey curve is given by $$N_7(E)= 2^{\alpha}\cdot3\cdot5\cdot7,\quad \alpha=0,1.$$
Hence $N=N_7(E)=105,210.$ Let $N=105$, the dimension of weight $2$ newforms of trivial Nebentypus character is $2$ and $c_{43}= 4, 4\sqrt{5}$. Since $43\mid XYZ$, the curve has a multiplicative reduction at $43$ and also $43$ is co-prime to $7N$, hence
$$\textrm{trace } \rho_{\ell}^{E}(Frob_p)= \pm (p+1).$$
Therefore the Fourier coefficient $c_p$ satisfies $7\mid \mathcal{N}_{K/\mathbb{Q}}(c_{43}\pm  44)$. A quick computation with possible choices of $c_{43}$ leads to a contradiction.

Now we consider $N=210$, the dimension of weight $2$ newforms of trivial Nebentypus character is $2$ and we find that $c_{43}= -4,8,-12$. This will lead to a contradiction except for the newform corresponding to $210.2.a.c$ in {\it The $L$-functions and Modular Forms Database} \cite{LMFDB}.

Since $29\mid XYZ$, using similar argument as above we conclude that $7\mid (c_{29}\pm  30).$ For the newform $210.2.a.c$ we have $c_{29}=6$ but $7\nmid (6\pm   30)$. This is a contradiction. 
\end{proof}


\section{\bf Preliminaries}

We consider \eqref{ES-oneterm} with $\ell$ prime. Write $x=\frac{n}{s}$ and 
$y=\frac{m}{t}, m\neq 0$ with $s,t$ positive integers and gcd$(n,s)=$ gcd$(m,t)=1.$ Then 
\eqref{ES-oneterm} becomes
$$(n+s)\cdots (n+(i-1)s)(n+(i+1)s)\cdots (n+ks)= \frac{s^{k-1}m^\ell}{t^\ell}.$$
Since the left hand side is an integer and gcd$(n, s)=\gcd(m, t)=1$,  we get $s^{k-1}=t^\ell.$
Since $\ell$ is prime, we have $\gcd(\ell, k-1)=\ell/\ell_0$ with $\ell_0\in \{1, \ell\}$. Hence there 
is a positive integer $d$ such that $s=d^{\ell_0}$ and $ t= d^{(k-1)\ell_0/\ell}$ and 
therefore \eqref{ES-oneterm} gives rise to the equation
\be\label{AP-oneterm}
(n+d^{\ell_0})\cdots (n+(i-1)d^{\ell_0})(n+(i+1)d^{\ell_0})\cdots (n+kd^{\ell_0})=m^{\ell}
\ee
with gcd$(n,d)=1.$ Thus the problem of finding rational points in \eqref{ES-oneterm} converts 
into finding integral solutions of the equation \eqref{AP-oneterm} which is 
equivalent to finding perfect powers in 
a product of consecutive terms of an arithmetic progression with one term missing.  When 
$d=1$, it is known by the results of Saradha and Shorey \cite{SaSh1} and \cite{SaSh3} that 
the only integral solutions of the equation \eqref{AP-oneterm} are given by 
$$\frac{6!}{5}=12^2, \quad \frac{10!}{7}=720^2 \quad {\rm and} \quad  \frac{4!}{3}=2^3.$$
Hence we consider $d>1$. Also we assume that $1<i<k$ since the case of $i\in \{1, k\}$ is completely 
solved in \cite{BBGH}. By symmetry, we further assume that $1<i\leq \frac{k+1}{2}$. 

We observe that any prime which divides at most one term in the product on the left hand 
side of \eqref{AP-oneterm} can occur only to an  $\ell$-th power. 
\textcolor{\defc}{If a prime $p$ divides both $n+jd^{\ell_0}$ and $n+rd^{\ell_0}$, where $r,j\neq i$, $r\neq j$, then $p$ divides $(j-r)d^{\ell_0}$. If $p$ divides $d$, then $p$ has to divide $n$, which is not possible as $\gcd(n,d)=1$. Hence $p$ must divide $(j-r)$, implying that $p<k$.}
Hence we can write each term as 
$$n+jd^{\ell_0}=a_jx_j^\ell, \quad 1\leq j\leq k, j\neq i$$
with $P(a_j)<k$ and every prime dividing $a_j$ occurs in another $a_r$ with
$r\neq i,j.$ From $\gcd(n, d)=1$, we observe that 
$$\gcd(a_j, d)=1, \quad  1\leq j\leq k, j\neq i.$$
Also for $\ell \geq 3$, $a_j$'s are taken as positive, by merging the negative sign in $x_j,$ 
if necessary. When $\ell=2,$ then we will consider $a_j$'s with necessary signs. Further, since $3\leq k\leq 8,$ we write 
$$a_j=2^{\alpha_j}3^{\beta_j}5^{\gamma_j}7^{\delta_j}$$
with integers $0\leq \alpha_j,\beta_j,\gamma_j,\delta_j<\ell$. Note that by \eqref{AP-oneterm},
$$\sum{\alpha_i}=\sum{\beta_i}=\sum{\gamma_i}=\sum{\delta_i}\equiv 0 \ \pmod{\ell}.$$
For a prime $2\leq p\leq k$ with $p\nmid d$, let $j_p$ be the least $j$ with $1\leq j\leq k, j\neq i$ such that 
 $$\nu_p(n+j_pd^{\ell_0})\geq \nu_p(n+jd^{\ell_0}) \quad {\rm for} \quad  1\leq j\leq k, j\neq i.$$ 
Then for $j\neq j_p$, we get from $n+jd^{\ell_0}-(n+j_pd^{\ell_0})=(j-j_p)d^{\ell_0}$ that 
$$\nu_p(n+jd^{\ell_0})\leq \nu_p(j-j_p) \quad {\rm for} \quad  1\leq j\leq k, j\neq i, j_p.$$ 
and equality holds when $\nu_p(n+jd^{\ell_0})\neq \nu_p(n+j_pd^{\ell_0})$. Further we observe that for $k/2<p<k$, either 
$p\nmid a_j$ for any $j$ or  there is a $1\leq j<j+p\leq k$ with either 
$$\left(\nu_p(a_j)=1, \quad \nu_p(a_{j+p})=\ell-1\right) \quad {\rm or} \quad \left(\nu_p(a_j)=\ell-1, \quad \nu_p(a_{j+p})=1\right).$$ 

As in \cite{BBGH}, we will form ternary equations of the shape
\begin{equation}\label{abcxyz}
ax^\ell+by^\ell+cz^{\ell}=0 \quad {\rm with} \quad 1\leq a\leq b\leq c \ {\rm and} \  \gcd(a, b, c)=1.
\end{equation} 
or of the shape
\begin{equation}\label{abc2}
ax^\ell+by^\ell=cz^2 \quad  {\rm with} \quad a, b, c\geq 1 \ {\rm and} \  \gcd(a, b, c)=1.
\end{equation}
This is achieved by considering the following four  identities. 
\begin{enumerate}[  ]
\item Let $1\leq p\neq q\leq k.$  Then 
$$a_qx^\ell_q-a_px^\ell_p=(n+qd^{\ell_0})-(n+pd^{\ell_0})=(q-p)d^{\ell_0}$$
gives rise to \eqref{abcxyz} when $\ell_0=\ell.$
\item Let  $p,q,r$ be three distinct integers  with $1\leq p,q,r\leq k.$ Then 
\begin{align*}
&(r-q)a_px^\ell_p+(p-r)a_qx^\ell_q+(q-p)a_rx^\ell_r\\
=&(r-q)(n+pd^{\ell_0})+(p-r)(n+qd^{\ell_0})+(q-p)(n+pd^{\ell_0})=0
\end{align*}
gives rise to \eqref{abcxyz}.
\item Let $p,q,r,t$ be four distinct integers with $1\leq p,q,r,t\leq k$ with $p+q=r+t.$ Then  
\begin{equation}\label{pqrt}
a_pa_q(x_px_q)^\ell-a_ra_s(x_rx_s)^\ell\\
\end{equation}
$$=(n+pd^{\ell_0})(n+qd^{\ell_0})-(n+rd^{\ell_0})(n+td^{\ell_0})=(pq-rt)d^{2\ell_0}$$
gives rise to \eqref{abc2}. When $\ell_0=\ell,$ then  \eqref{pqrt} gives rise to 
\eqref{abcxyz}.
\end{enumerate}
We shall denote these equations by $[p,q]$, $[p,q,r]$,  $[p,q,r,t].$ It will be clear from the usage whether we consider $[p,q,r,t]$ as having signature $(\ell,\ell,\ell)$ or $(\ell,\ell,2).$

\section{\bf $k\geq 4,\ell>2$ and $\gcd(\ell,k-1)=1$}\label{gcd1}

Throughout this section, we shall assume that  $4\leq k\leq 8, \ell>2$ and $\gcd(\ell, k-1)=1$ so that $\ell_0=\ell$. \textcolor{\defc}{With these assumptions we rewrite the 
equation \eqref{AP-oneterm} as} 
\begin{align}\label{d-ell}
(n+d^{\ell})\cdots (n+(i-1)d^{\ell})(n+(i+1)d^{\ell})\cdots (n+kd^{\ell})=m^{\ell}
\end{align}
with $\gcd(n, d)=1$ and $1<i\leq \frac{k+1}{2}$ so that \eqref{d-ell} does not hold thereby proving Theorem \ref{Th2} for $\ell\geq 3.$
Our strategy is to form a ternary equation as in \eqref{abcxyz} or \eqref{abc2} which will yield no solution. The following lemma is very useful.
\begin{lemma}\label{I&II}
Let $1\leq j\le k, j\neq i$ be satisfies one of the following:
\begin{enumerate}[(i)]
\item $P(a_ja_{j+2^s})\leq 2$ for some $s\geq 0$ 
\item $P(a_ja_{j+3^s})\leq 3, 2\nmid a_ja_{j+3^s}$ for some $s=0,1$
\item $P(a_ja_{j+3t})=3, \nu_3(a_j)=\nu_3(a_{j+3t})=1$ with $t\in \{1, 2\}$
\item $P(a_ja_{j+t}a_{j+2t})\leq 3$ for some positive integer $t=2^s$ with $s\geq 0$
\item $P(a_ja_{j+t}a_{j+2t})=3$ for some positive integer $t=3\cdot 2^s$ with $s\geq 0$.
\end{enumerate}
Then \eqref{d-ell} has no solution. 
\end{lemma}

\begin{proof}
In case of $(i)$, we form $[j, j+2^s]$ and this has no solution by 
Lemma \ref{eq1} with $p=2$. In case of $(ii)$, we form $[j, j+3^s]$ and this 
has no solution by Lemma \ref{eq1} with $p=3$. In case of $(iii)$, we form 
\emph{    }$[j, j+3t]$ and the resulting equation $ax^\ell+by^\ell+cz^\ell=0$ satisfy 
$P(abc)\leq 2$ and this has no solution by Lemma \ref{eq1} with $p=2$. In the case $(iv)$, we have either 
$$3\nmid a_ja_{j+t} \quad {\rm or} \quad 3\nmid a_ja_{j+2t} \quad {\rm or} \quad  3\nmid a_{j+t}a_{j+2t}.$$ 
Any of this gives rise to a case of the form $(i)$ since $P(a_ja_{j+t}a_{j+2t})\leq 3$ and hence 
there is no solution by Lemma \ref{eq1} with $p=2$. Consider the case $(v)$. Since 
$P(a_ja_{j+t}a_{j+2t})=3$, we have $3\nmid d$ and each of $n+jd^\ell$, $n+(j+t)d^\ell$ and $n+(j+2t)d^\ell$ 
is divisible by $3$. Since $t=3\cdot 2^s$, we see that 9 divides at most one of these three terms.
Let $\ve_1<\ve_2$ with $\varepsilon_1, \ve_2\in \{0, 1, 2\}$ be such that  
$3||(n+(j+\ve_1t)d^\ell)$ and $3||(n+(j+\ve_2t)d^\ell)$. Then 
$$3||a_{j+\ve_1t}\ {\rm and}\ 3||a_{j+\ve_2t}.$$ 
Now \emph{    }$[j+\ve_1t, j+\ve_2t],$  after cancellation of common factors gives rise to an equation as in Lemma \ref{eq1} with $p=2$ and hence has no solution.
\end{proof}
\begin{lemma}
Equation \eqref{d-ell} has no solution whenever $d$ is even.
\end{lemma}
\begin{proof}
Suppose $2|d$. Then $2\nmid a_j$ for any $j$ with $1 \leq j \leq k, j\neq i.$ Let $k=4$. We have 
$P(a_1a_4)\leq 3$ and  Lemma \ref{I&II} (ii) with $(j,s)=(1,1)$ implies that \eqref{d-ell} has no solution. 
For $k\in \{5, 6\}$ and $k\in \{7, 8\}$ with $i\neq 4$, we get $p\nmid a_4a_5$ for 
$p\in \{5, 7\}$. Since  $2\nmid a_4a_5,$ Lemma \ref{I&II} (ii)  with $(j,s)=(4,0)$ implies that \eqref{d-ell} has no solution.  
Let $k\in \{7, 8\}$ and $i=4$. We have  $P(a_2a_3)\leq 3$ and $P(a_5a_6)\leq 3$ according as 
$5\nmid a_2a_3$ and $5|a_2a_3$, respectively. Now the assertion follows by Lemma \ref{I&II} (ii) with 
$(j, s)=(2, 0), (5, 0)$, respectively. 
\end{proof}
From now on we shall assume that $2\nmid d$. For $k\geq 5,$ observe that 2 divides at least two $a_i$' s since $2\nmid d.$  

\begin{lemma}\label{  4}
Let the assumptions of this section hold. \\
(i) Equation [2,7,1,8] implies that $\ell \leq 5$ whenever
$$5|a_2a_7\ or\ 7|a_1a_8\ and \ 5\nmid a_1a_8.$$
If also 
$$ (a) \ 3\nmid a_1a_2a_7a_8,\ then \ 4|a_2a_7\ or \ 4|a_1a_8$$
$$ (b)\ 3|a_2a_5a_8, then \ (\beta_2,\beta_8)\neq (1,1).$$ 
(ii)  Equation  [2,7,1,8]  implies that $\ell \leq 7$ if $35|a_1a_8.$
\end{lemma}

\begin{proof}
Note that
\begin{equation}\label{22718}
2||\gcd(a_2a_7,a_1a_8).
\end{equation}
Hence $[2,7,1,8]$ reduces to 
$$Ax^\ell+By^\ell=3d^{2\ell}\ with\ P(AB)\leq 3$$
and
$$5|xy\ \textrm{ if } 5|a_2a_7\ \textrm{ or } 7|xy\ \textrm{ if }\ 7|a_1a_8\ \textrm{ and }\ 5\nmid a_1a_8.$$
By Lemma \ref{eq4}(iv) we get that $\ell\leq 5$ except perhaps when $P(AB)=2$ which is the case when $3\nmid a_1a_2a_7a_8.$ In this case by \eqref{22718} the equation reduces to 
$$x^\ell+2^\alpha y^\ell=3d^{2\ell}, 5|xy, \alpha\geq 0$$
which by Lemma \ref{eq4}(ii) gives $\ell\leq 5.$
This proves the first assertion of (i). When $\alpha=0,$ by Lemma \ref{eq1}, equation \eqref{d-ell} has no solution. Thus we may assume that $\alpha\geq 1$ which implies that $4|a_2a_7$ or $4|a_1a_8.$
This proves (a). 
Suppose $3|a_2a_5a_8$ and $(\beta_2,\beta_8)=(1,1).$  Then  [2,7,1,8] reduces to $x^\ell+2^\alpha y^\ell=z^\ell$
which is excluded by Lemma \ref{eq1} and this proves (b).

(ii) Here $[2,7,1,8]$ reduces to 
$Ax^\ell+By^\ell=3d^{2\ell}\ with\ P(AB)\leq 5, 7|xy$
and hence the result follows from \ref{eq4} (v).
\end{proof}

We consider different values of $k$.

\subsection{\bf Let $k=4$.}  Then $i=2, P(a_j)\leq 3$ and $\ell\geq 5$ since $\gcd(\ell, 3)=1$.   
Also 
$$ \ \alpha_1+\alpha_3, \  \beta_1+\beta_4\equiv 0 \pmod{\ell}. 
$$ 
Hence we have the following cases.
\begin{center}
\begin{tabular}{|c|l|c|c| c||}
\hline
$(i)$ & $2|(n+4d^\ell), 3\nmid (n+4d^\ell)$ & $a_4=1, 2\nmid a_3$ &  $a_3a_4=3^s, s\geq 0$\\
\hline
$(ii)$ & $2|(n+4d^\ell), 3|(n+4d^\ell)$ & $a_3=1, 2\nmid a_4$ & $a_3a_4=3^s, s\geq 0$ \\
\hline
$(iii)$ & $2\nmid (n+4d^\ell), \al_1=\ell-1$ & $a_4=3^\beta$ &  $[1, 4]$\\
\hline
$(iv)$ & $2\nmid (n+4d^\ell), \al_3=\ell-1$ & $a_4=3^\beta$ & $[3, 4]$\\  \hline 
\hline
\end{tabular}
\end{center}
The cases $(i)$ and $(ii)$ are excluded by Lemma \ref{I&II} $(ii)$.  The cases $(iii)$ and $(iv)$ are excluded by 
Lemma \ref{eq2} with $(a,b,c)=(-a_1,a_4,3)$ or $(-a_3,a_4,1)$ since $\ell\geq 5$.

\subsection{\bf Let $k=5$.}  Then $\ell\geq 3$ and $i\in\{2,3\}$. Since $5|(n+jd^\ell)$ for at most one $j$, for each $j\neq i,$ we have $\gamma_j=0$ and hence $P(a_j)\leq 3.$  If $i=2$, then Lemma \ref{I&II} $(iv)$ with $j=3, t=1$ gives the assertion. Hence we suppose that 
$$i=3.$$ 
By Lemma \ref{GHP}, we have 
$$\ell\leq 5.$$ By considering Lemma 
\ref{I&II} $(i)$ with $j=4, s=0$, we may assume that 
$$3|a_4a_5.$$ 
We have 
$$\al_1+\al_5\equiv \al_2+\al_4\equiv 0 \pmod{\ell}.$$
Equation $[1, 5, 2, 4]$ gives rise to \eqref{abcxyz} with $abc=3^s$, $s= 0, 1$ which is excluded by Lemma \eqref{eq1} with $p=3$.

\subsection{\bf Let $k=6$.} Then $\ell\geq 3, \ell\neq 5$ and $i\in \{2, 3\}$. Also  $P(a_j)\leq 5$.  If $5\nmid a_j$ for any $j$, then 
Lemma \ref{I&II} $(iv)$ with $j=4, t=1$ gives the assertion. Thus $5|a_j$ for some $j$ implying 
$$5|a_1,a_6$$ 
with 
$(\gamma_1, \gamma_6)\in \{(1, \ell-1), (\ell-1, 1)\}$ and $5\nmid a_j$ for $j\neq 1, 6$. If $i=2$, then 
Lemma \ref{I&II} $(iv)$ with $j=3, t=1$ gives the assertion. So we suppose that $$i=3.$$ 
Further from Lemma 
\ref{I&II} $(i)$ with $j=4, s=0$, we may assume that 
$$3|a_4a_5.$$ 
Suppose $3|a_5$. Then $\nu_3(a_2a_5)=\ell$ and $3\nmid a_1a_6$. Considering $[2, 5, 1, 6]$, we have from 
\begin{equation}\label{2516}
a_2a_5(x_2x_5)^\ell-a_1a_6(x_1x_6)^{\ell}=4d^{2\ell}=4(d^2)^\ell, 
\end{equation}
an equation $ax^\ell+by^\ell=(d^2)^\ell=z^\ell$ with $P(ab)\leq 2$. This is excluded by Lemma \ref{eq1} with $p=2$.  
Hence 
$$3|a_4.$$
Then $3\nmid a_2a_5a_6$ and $(\beta_1, \beta_4)\in \{(1, \ell-1), (\ell-1, 1)\}$.  Also note that $\al_2+\al_4+\al_6\equiv \al_1+\al_5\equiv 0 \pmod{3}$. 

\textcolor{\defc}{We first let $\ell=3$}. In Table 1 we list the  cases depending on different choices of $j_2$
and the choices of $p,q$ to form the equation $[p,q].$ 
The equations in the last column of Table 1 are excluded by Lemma \ref{ell3}. 
\begin{table}
\begin{center}
\begin{tabular}{|c|c|c|c|}
\hline
$j_2$ & Cases & $[p,q]$ & Equation\\
\hline
1 & $a_4=9, a_2=1$  &  $[2, 4]$  &  $x^3+2y^3+9z^3=0$ \\ \hline
1 & $a_4=3, a_2=1$ &  $[1, 2]$  &  $x^3+y^3+cz^3=0$ \\
& $a_1\in \{3^2\cdot 5, 3^2\cdot 5^2\}$ & & $c\in \{45, 225\}$ \\ \hline
2 & $a_2=a_5=1$ & $[2, 5]$ &  $x^3+y^3+3z^3=0$ \\ \hline
4 & $a_2=2, a_4=6$  &  $[2, 4]$  &  $x^3+y^3+3z^3=0$ \\ \hline
4 & $a_2=2, a_4=18$  &  $[4, 5]$  &  $x^3+y^3+18z^3=0$ \\ \hline
5 & $a_6=5, a_5=2$ &  $[5, 6]$ &  $x^3+2y^3+5z^3=0$ \\  \hline
5 & $a_6=5^2, a_2=1$ & $[2, 6]$ &  $x^3+4y^3+25z^3=0$ \\ \hline
6 & $a_5=1, a_6\in \{5, 5^2\}$  & $[5, 6]$  &  $x^3+y^3+cz^3=0$ \\  
&&&  $c \in \{5, 25\}$ \\ \hline 
\end{tabular}
\caption{\label{tab1:table-name} {$(k,\ell,i)=(6,3,3)$}}
\end{center}
\end{table}
\textcolor{\defc}{Thus we can consider}
$$\ell\geq 7.$$
 We have 
$\al_{j_2}\geq \ell-3$. Let $j_2=4$. Then $a_5=1$, $P(a_4)\leq 3$ and $\al_4\geq \ell-3\geq 4$. This gives 
equation $[4, 5]$ in  the form $x^\ell+y^\ell=a_4z^\ell$ and this is not possible by Lemma \ref{eq2}. Let $j_2\in \{1, 2, 5, 6\}$. We consider $[2, 5, 1, 6].$ We 
obtain from \eqref{2516} the equation $ax^\ell+by^\ell=(d^2)^\ell=z^\ell$ with $P(ab)\leq 3$ and 
$\nu_2(ab)\geq (\ell-3)-2$. If $\ell\geq 11$, we get a contradiction by Lemma \ref{eq2}. Thus 
$$\ell=7.$$  
We rewrite the equation as $ax^\ell+by^\ell=(d^\ell)^2=z^2$ with 
$P(ab)\leq 3$. Also  $\nu_2(ab)\leq( \ell-2)-2=3$. By Lemma \ref{eq4} $(iv)$, we get a contradiction. 

 {\subsection{\bf Let $k=7$.} 

We have  $\ell\geq 5$ and $i\in \{2, 3, 4\}$. Also  $P(a_j)\leq 7$.  Since $7|(n+jd^\ell)$ for at most one $j$, we have 
$\delta_j=0$ for each $j$. Hence $P(a_j)\leq 5$ for each $j$. If 
$5\nmid a_j$ for any $j$,  then Lemma \ref{I&II} $(iv)$ with $j=5, t=1$ gives the assertion. 
Thus $5|a_j$ for some $j$ implying 
$$5|a_1, 5|a_6 \  {\rm with} \  (\gamma_1, \gamma_6)\in \{(1, \ell-1), (\ell-1, 1)\} \quad {\rm and} \quad  
5\nmid a_j \ {\rm for} \ j\neq 1, 6$$ 
or  
$$5|a_2, 5|a_7 \  {\rm with} \  (\gamma_2, \gamma_7)\in \{(1, \ell-1), (\ell-1, 1)\} \quad {\rm and} \quad  
5\nmid a_j \ {\rm for} \ j\neq 2, 7.$$

First we take $5|a_2$ and $5|a_7$. Then $i=3,4$ and Lemma \ref{I&II} $(iv)$ with $j=4, t=1$ when $i=3$ and 
with $j=1, t=2$ when $i=4$ shows that there is no solution. 

Secondly we take  $5|a_1$ and $5|a_6$. If $i\in \{2, 4\},$ 
then Lemma \ref{I&II} $(iv)$ with $j=3, t=1$ when $i=2$ and with $j=3, t=2$ when $i=4$ shows that there is no solution. Therefore we suppose that 
$$i=3.$$  
Further from Lemma \ref{I&II} $(i)$ with $j=4, s=0$, we may assume that 
$$3|a_4a_5.$$ 
Let $3|a_5$. Then $\nu_3(a_2a_5)=\ell$ and $3\nmid a_1a_6$. Considering 
$[2, 5, 1, 6]$, we have  \eqref{abcxyz} with $abc=2^s$, $s\geq 0$ which is excluded by Lemma \eqref{eq1} with 
$p=2$.  Thus we assume that 
$$3|a_4.$$ 
Then $3|a_1a_4a_7$ and $3\nmid a_2a_5a_6$. 
If $\beta_1=\beta_7=1$, then Lemma \ref{I&II} $(iii)$ with $j=1, t=2$ gives a contradiction. Hence 
$$(\beta_1, \beta_4, \beta_7)\in \{(1, 1, \ell-2), (\ell-2,1,1)\}.$$ 

Let $j_2\in \{1, 2, 5, 6\}$. Then $\al_{j_2}\geq \ell-3$. Then $[2, 5, 1, 6]$ 
gives the necessary contradiction as in the case $k=6$,  $\ell> 7$. Thus we need to consider
$$\ell\in\{5,7\}.$$
First let $\ell=5$. Then 
$$\al_{j_2}=\al_{j_2+4}=2\ \text{if} \ j_2\in \{1, 2\}$$
and  
$$\al_{j_2}=\al_{j_2-4}=2\ \text{if} \ j_2\in \{5, 6\}.$$
Also 
$$\beta_1=1 \ \text{or} \ 3.$$ From  \eqref{2516}, we get 
an equation of the form $x^5+y^5=3^\alpha z^5, \alpha=1,3$ which is not possible by Lemma \ref{eq1} with $p=3$. 
Next let $\ell=7$. Then from $\al_{j_2}\geq 7-3=4, \beta_1=1, 3\nmid a_2a_5a_6$ and equation 
\eqref{2516}, we obtain equations as in  Lemma \ref{eq4} $(iv)$ which have no solution.  

Let $j_2\in \{4, 7\}$. Then $\al_{j_2}\geq \ell-2$.  We consider $[4, 5, 2, 7]$ which is 
an equation of the form $ax^\ell+by^\ell=(d^2)^\ell$ with $P(ab)\leq 3$ and $\nu(ab)\geq \ell-3.$
By Lemma \ref{eq2}, we may therefore suppose that 
$$\ell=5.$$  
Then $a_5=1$ if $j_2=4$ and $a_2=1$ if 
$j_2=7$. Since $P(a_1a_4)\leq 5$, we consider $[4, 5]$ if $j_2=4$ and 
$[1, 2]$ if $j_2=7$. We obtain equations of the form $x^5+y^5=a_4z^5$ and 
$x^5+y^5=a_1z^5$, respectively which are excluded by Lemma \ref{ell5}. This completes the case for  
$k=7$. 

{\subsection{\bf Let $k=8$.} 
\vskip 2mm
We have 
$$\ell\geq 3, \ell\neq 7, i\in \{2, 3, 4\}\ \textrm{ and } \ P(a_j)\leq 7.$$ 
Suppose $P(a_j)\leq 3$ for each $j$. Then Lemma \ref{I&II} $(iv)$ with $j=5, t=1$ 
gives the assertion. Thus either $5|a_j$ for some $j$ or $7|a_j$ for some $j.$
Further if 5 does not divide any $a_j,$ then at least two of $a_5,a_6,a_7$ have $P(a_i)\leq 2$ and this is excluded by Lemma \ref{I&II} $(i).$ Hence we may assume that
$$5|a_j, 5|a_{j+5} \ {\rm for \ some} \ j\in \{1, 2, 3\} \ {\rm and} \ (\gamma_j, \gamma_{j+5})\in \{(1, \ell-1), (\ell-1, 1)\}.$$

Let 
$$5|a_3a_8.$$
Then $i\in\{2,4\}$ and $P(a_5a_6a_7)\leq 3.$ Hence the assertion follows by Lemma \ref{I&II} $(iv)$ with $j=5, t=1$. 

Next let 
$$5|a_2a_7.$$ 
Then  $i\in \{3,4\}.$  Let $i=3.$ Then $P(a_4a_5a_6)\leq 3$ and the assertion follows by 
Lemma \ref{I&II} $(iv)$ with $ j=4,t=1.$ Next let $i=4.$ Then by Lemma \ref{  4}(i) we get 
$$\ell\leq 5.$$
Note that by Lemma \ref{I&II}(i), we may assume that
$3|a_5a_6.$
Hence either
$$ 3|a_3a_6\ \textrm{ or }\ 3|a_2a_5a_8.$$
Suppose 
$$3|a_3a_6.$$ 
Then 
$3\nmid a_1a_2a_7a_8.$ Hence by Lemma \ref{  4}(i)
we have 
\begin{equation}\label{4a2}
4|a_2a_7\ \textrm{ or }\ 4|a_1a_8.
\end{equation} 
If $a_5\in \{1, 2\}$, then we consider  $[5, 6]$ when $a_5=1$ and 
$[3, 5]$ when $a_5=2$ and  this is excluded by Lemma \ref{ell5} since $P(a_3a_6)=3$.  
Thus $4|a_5$ which together with \eqref{4a2} gives $4|a_1$ and $2||a_3, 2||a_7$ so that 
$$\al_3=\al_7=1, \al_1\geq 2\ \textrm{ and } \ \al_5\geq 2.$$
When $\ell=5,$ we have
$\al_1+\al_3+\al_5+\al_7\equiv 0 \pmod{5}$ giving $\al_1+\al_5\in \{3, 8\}$ which together with $\al_j<5$ implies that
$$\al_1=\al_5=4.$$
This means 
$$2^4|(n+5d^\ell-n-d^\ell)=4d^\ell$$
which is not possible and hence $\ell=5$ is excluded. Let $\ell=3$. 
We have $a_5=2^{\al_5}$ and $a_7=2^{\al_7}5^{\gamma_7}$ with $\gamma_7\in \{1, 2\}$ and either 
$\al_5=\al_7=0$ or $1\leq \al_5, \al_7\leq 2$. This is excluded by considering equation
$[5, 7]$ and Lemma \ref{ell3}. 

Next let 
$$3|a_2a_5a_8.$$ 
so that 
$$a_3=1\ \textrm{ or } \ a_6=1.$$ 
When $\ell=3,$ we have $\beta_2=\beta_8=1$ which is excluded by Lemma \ref{  4}(i)(b).
Let $\ell=5.$ We apply Lemma \ref{ell5} to equations 
$$[2, 3]\ \textrm{ when } \ a_3=1$$
 and 
$$[6, 7]\ \textrm{ when }  \ a_6=1$$ 
to get the assertion since  $P(a_2a_7)=5$.  This concludes the case $5|a_2a_7.$

Lastly we take 
$$5|a_1a_6.$$
Then 
$$i\in\{2,3,4\}\ \textrm{ and } \ P(a_j)\leq 3\ \textrm{ for }  \ j\notin \{1, 6, 8\}.$$ 
The assertion of the theorem follows 
$${\rm if }\ i=2\ {\rm by \ Lemma \ \ref{I&II}} (iv)\ {\rm with}\ j=3, t=1$$
and
$${\rm if}\ i=4\ {\rm by\ Lemma \ \ref{I&II}} (iv) \ {\rm with}\ j=3,t=2.$$ 
So we consider
$$i=3.$$
If $3\nmid a_4a_5$, then again the assertion follows by 
Lemma \ref{I&II} $(i)$ with $j=4, s=0$. Hence we assume that 
\begin{equation}\label{3a1}
3|a_1a_4a_7\ \textrm{ or } \ 3|a_2a_5a_8.
\end{equation}
 Let $\ell=3$. Then 
$$\beta_j=1\ \textrm{ for }  \ j\in \{1, 4, 7\}\  \textrm{ if }  \ 3|a_1a_4a_7$$
and 
$$\beta_j=1\ \textrm{ for }  \ j\in \{2, 5, 8\}\ \textrm{ if }  \ 3|a_2a_5a_8.$$
Also note that
\begin{equation}\label{a4a7}
P(a_4a_7a_2a_5)\leq 3.
\end{equation}
Thus $\ell=3$ is excluded 
$$\textrm{ if }  \ 3|a_1a_4a_7\ {\rm by\ Lemma\ \ref{I&II}} (iii)\ \textrm{ with } \ j=4, t=1$$ 
and 
$$ \textrm{ if } \ 3| a_2a_5a_8\ {\rm by\ Lemma\ \ref{I&II}}(iii)\ \textrm{ with }  \ j=2, t=1.$$ 
Thus we have
$$\ell\geq 5.$$ 
Let us now assume that  
$$7\nmid a_j\ \textrm{ for any }  \ j.$$
Then Lemma \ref{GHP} with 
$m=n+4d^\ell, s=d^\ell$ and \eqref{a4a7} gives $$\ell\leq 5.$$
If $3|a_2a_5a_8$,  then 
Lemma \ref{I&II} $(v)$ with $j=2, t=3$ gives the assertion. If $3|a_1a_4a_7,$
then
$$P(a_2a_5a_8)\leq 2$$ 
and either
 $$a_2=1\ \textrm{ or } \  a_5=1\ \textrm{ or } \ a_8=1.$$
The case $\ell=5$ is excluded by taking
 $$[1, 2]\ \textrm{ if } \ a_2=1;  [5, 6]\ \textrm{ if } \ a_5=1;[7, 8]\ \textrm{ if } \ a_8=1.$$
Hence we may assume that
\begin{equation}\label{7a1}
7|a_1a_8.
\end{equation}
Thus $35|a_1a_8.$ So by Lemma \ref{  4}(i) we have $\ell=5$ since $\ell\neq 7.$
Then by \eqref{3a1}, \eqref{7a1} and Lemma \ref{ell5},$[2, 7, 1, 8]$ reduces to  
$$
a x^5+(d^2)^5=b(-y)^5
$$
with $a=a_2a_7/6>1.$
Hence we assume that
$${\rm either} \qquad 2\bigg|\frac{a_2a_7}{6} \qquad {\rm or} \qquad 3\bigg|\frac{a_2a_7}{6}.$$ 
Note that $\gcd(a_2, a_7)=1$. In Table 2 we give different possibilities of $a_2$ and use equation $[p, q]$  with suitable values of $p,q$ to exclude these possibilities by Lemma \ref{ell5}.
\begin{table}
\begin{center}
\begin{tabular}{|c|c|c|}
\hline
$a_2$ & Cases  & $[p,q]$\\
\hline
$2\mid a_2 ,3\nmid a_2$ & $a_4=2^{\alpha_4}3^{\beta_4},\beta_4>0, a_5=1$& $[4,5]$\\
\hline
$2\mid a_2, 3|a_2$ &$a_6=2^{\alpha_6}5^{\gamma_6},\gamma_6>0, a_7=1$ &$ [6,7]$\\
\hline
$2\nmid a_2,3|a_2 $ & $a_4=1, a_5=2^{\alpha_5}3^{\beta_5},\beta_5>0$ & $[4,5]$\\
\hline
$2\nmid a_2, 3\nmid a_2 $ &$a_1=2^{\alpha_1}3^{\beta_1}5^{\gamma_1}7^{\delta_1},\gamma_1>0,\delta_1>0, a_2=1$ & $[1,2]$\\
\hline
\end{tabular}
\caption{\label{tab2:table-name} {$(k,\ell,i)=(8,5,3)$}}
\end{center}
\end{table}

This concludes the case $k=8.$

\section{\bf $\ell>2$ and $\gcd(k-1, \ell)>1$}\label{gcd>1}

In this section, we consider $\ell>2$ and $\gcd(\ell, k-1)>1$ so that 
$$(k, \ell)\in \{(4, 3), (6, 5), (7, 3), (8, 7)\}.$$
Further from $\ell_0=1$, equation \eqref{AP-oneterm} becomes 
\begin{align}\label{d-1}
(n+d)\cdots (n+(i-1)d)(n+(i+1)d)\cdots (n+kd)=m^{\ell}
\end{align}
with $\gcd(n, d)=1$ and $1<i\leq \frac{k+1}{2}$. 
Recall that for $p\nmid d$ with $1\leq p\leq k$, $j_p$ is the least $j\neq i$ such that 
$\nu_p(n+jd)\leq \nu_p(n+j_pd)$ for all $j\neq i$. 

\begin{lemma}\label{l^2}
Given $\ell$ prime $>2$, let $q$ be either a prime or a power of $\ell$ so that 
$\ell|\varphi(q)$. Let $\Lambda_q$ be the set of solutions of 
$x^{\frac{\varphi(q)}{\ell}}\equiv 1 \pmod{q}$. 
When $q|d$, we have for $j_1\neq j_2$, 
\begin{align}\label{q|d}
a_{j_1}\equiv \la_qa_{j_2}  \pmod{q} \quad  for \ some \quad  \la_q\in \Lambda_q.
\end{align}
Let $0\leq j_q<q$ be such that $q|(n+j_qd)$. For 
$1\leq j_1\neq j_2\leq k$ and $j_1\neq j_q, j_2\neq j_q$, we have  
\begin{align}\label{chi}
\frac{a_{j_1}}{j_1-j_q}\equiv \frac{\la_qa_{j_2}}{j_2-j_q} \pmod{q}
\  for \ some \quad  \la_q\in \Lambda_q.
\end{align}
\end{lemma}

\begin{proof}
When $q|d$, we have $a_jx^\ell_j=n+jd\equiv n \pmod q$ and hence $(a_j)^{\frac{\varphi(q)}{\ell}}\equiv 
n^{\frac{\varphi(q)}{\ell}} \pmod q$ and the assertion \eqref{q|d} follows. Let $q\nmid d$. From 
$a_jx^\ell_j=n+jd=n+j_qd+(j-j_q)d$, the assertion\eqref{chi} follows by observing 
$$(a_j)^{\frac{\varphi(q)}{\ell}}\textcolor{\defc}{\equiv}(a_jx^\ell_j)^{\frac{\varphi(q)}{\ell}}\equiv 
(j-j_q)^{\frac{\varphi(q)}{\ell}}d^{\frac{\varphi(q)}{\ell}} \pmod{q}.$$
\end{proof}
Here are some values of $\ell, q$ and $\Lambda_q$. 
\begin{align}\label{Laq}
\Lambda_q=\begin{cases}
\{ \pm 1, \pm  7\} & {\rm if} \ \ell=5, q=5^2;\\
\{ \pm 1, \pm  18,\pm   19\} & {\rm if} \ \ell=7, q=7^2;\\
\{ \pm 1, \pm  12\} & {\rm if} \ \ell=7, q=29.
\end{cases}
\end{align}

\begin{lemma}\label{l^21}
Let $\ell\in \{5, 7\}$ and $\ell^3|(n+j_\ell d)$. Let $\ve_\ell$ be given by $a_{j_\ell+\ell}=\ell\ve_\ell$ 
or $a_{j_\ell-\ell}=\ell\ve_\ell$ according as $1\leq j_\ell<j_\ell+\ell\leq k$ or 
$1\leq j_\ell-\ell<j_\ell\leq k$, respectively. Then 
for $1<j<k, j\neq j_\ell, j_\ell \pm \ell$, we have 
\begin{align}\label{epl}
\left(\frac{a_j}{j-j_\ell}\right)^{\ell-1}\equiv (\ve_\ell)^{\ell-1}  \pmod{\ell^2}.
\end{align}
\end{lemma}

\begin{proof}
Since $\ell^3|(n+j_\ell d)$, we have $\nu_\ell(n+(j_\ell \pm  \ell) d)=1$ and hence 
$\ell\nmid \ve_\ell$ and $(\ve_\ell)^{\ell-1}\equiv d^{\ell-1}  \pmod{\ell^2}$.  
Also for 
$1<j<k, j\neq j_\ell, j_\ell \pm \ell$, we have 
$$a_jx^\ell_j=n+jd=(n+j_\ell d)+(j-j_\ell)d\equiv (j-j_\ell)d  \pmod {\ell^2}.$$
Now the assertion follows by taking $(\ell-1)-th$ powers on both sides. 
\end{proof}

We now consider different cases of $(k, \ell)$. Recall that $1<i\leq \frac{k+1}{2}$.

\section{\bf The case $(k, \ell)=(4, 3)$}
 
 Let $(k, \ell)=(4, 3)$. Then $i=2,3.$ First let $i=2.$ 
We consider $[1, 3, 4]$ which gives $a_1x_1^3+2a_4x^3_4=3a_3x^3_3$ with $P(a_1a_2a_3)\leq 3$.\\
$\text{Let}\ D_{1,3,4}: \{(X,Y,Z): a_1X^3+3a_3Y^3+2a_4Z^3=0 \}  \ \text{and}\\
E_{6a_1a_3a_4}: x^3+y^2z+6a_1a_3a_4yz^2=0.$\\

We consider the morphism
\begin{align*}
 &D_{(1,3,4)} \xrightarrow{} E_{6a_1a_3a_4}\\
&(X,Y,Z) \xrightarrow{} (a_1^3(3a_3)^4x^4,a_1^3 (3a_3)^4y^2z,a_1^2 (3a_3)^4yz^2)
\end{align*}
This leads to  
$$a_1^{10}{(3a_3)}^{12}x^{12}+ a_1^{9}{(3a_3)}^{13}y^6z^3+2a_1^6{(3a_3)}^{12}a_4y^3z^6=0$$
which gives
$$a_1^{9}(3a_3)^{12}x^{12}+ a_1^{8}(3a_3)^{13}y^6z^3+da_1^4(3a_3)^{11}y^3z^6=0$$
where $d=6a_1a_3a_4.$
Now we substitute, 
$$x=a_1^3(3a_3)^4x^4, y=a_1^4(3a_3)^5y^3, z= (3a_3)^3z^3.$$
This gives $E_d: x^3+y^2z+dyz^2=0.$  We also note that $E_d$ and $E_{d_1}$ are isomorphic if $d/{d_1}$ is a cube. Hence it is enough to consider the case $d=6$ and $a_1=a_3=a_4=1.$
That is we need to solve
the  equation
\begin{equation*}
    (n+d)(n+3d)(n+4d)=y^3
\end{equation*}
with $(n+d)=y_1^3, (n+3d)=y_3^3, (n+4d)=y_4^3$, and $ y_1, y_3, y_4$ pairwise co-prime. This gives the equation 
\begin{equation}\label{eq-4-3-2}
y_1^3 + 2y_4^3 = 3y_3^3
\end{equation}
$y_1, y_3, y_4$ pairwise co-prime. The curve \eqref{eq-4-3-2} bi-rationally equivalent to the Weierstrass curve $V^2=U^3- 432\times 6^2$ and it has rank $1$. Hence, we conclude that the if the equation $$(n+d)(n+3d)(n+4d)=y^3$$ has a solution then it would arise from \eqref{eq-4-3-2}. 
Now, let $i=3$. Hence we have
$$(n+d)(n+2d)(n+4d)=y^3.$$
Using the substitution, $N=-(n+5d), Y=-y$, we have $$(N+d)(N+3d)(N+4d)=Y^3.$$
Since, $(n,d)=1$, we have $(N,d)=1$, therefore it is enough to consider the case for $i=2.$

\section{\bf The case $(k, \ell)=(6, 5)$}  

\subsection{\bf Let  $(k,\ell)=(6,5)$ and $i=2$.}  We have either 
$\al_j=0$ for all $j\neq 2$ or 
$$
\begin{matrix}
\al_1 & \al_3 & \al_5 & \al_4+\al_6 \\
1 & 3 & 1 & 0 \\
2 & 1 & 2  & 0
\end{matrix}
\quad {\rm or} \quad 
\begin{matrix}
\al_4 & \al_6 & \al_1+\al_3+\al_5 \\
4 & 1 & 0 \\
1 & 4 & 0
\end{matrix} 
$$
and either $\beta_j=0$ for all $j\neq 2$ or 
$$
\begin{matrix}
\beta_1 & \beta_4 & \sum_{j\neq 1, 4}\beta_j \\
1 & 4 & 0 \\
4 & 1 & 0 
\end{matrix}
\quad {\rm or} \quad 
\begin{matrix}
\beta_3 & \beta_6 & \sum_{j\neq 3, 6}\beta_j\\
1 & 4 & 0 \\
4 & 1 & 0.
\end{matrix} 
$$
Observe that 
$$\beta_5=0.$$
 Also 
$$5\nmid a_j\ \text{for} \ 1<j<6, j\neq 2.$$
 By Lemma \ref{k5}, we have either 
$$3|a_1a_4\ \text{or}\ 3|a_3a_6.$$
Let $2\nmid a_j$ for all $j$. Then $a_1=a_5=1$ if $3|a_3a_6$ and $a_3=a_5=1$ if $3|a_1a_4$. So 
$$[1, 3, 5]\  \text{if}\ 3|a_3a_6$$ 
or
$$[3,4, 5]\  \text{if}\ 3|a_1a_4$$ 
give equations of the form $x^5+y^5=cz^5$ with $P(c)\leq 3$ and are excluded 
by Lemma \ref{ell5}. Thus 
$$2|a_j\ \text{for some} \ j\neq 2.$$
Now we suppose that, $6|a_4$, then $a_3=a_5=1$ and $[3,4,5]$ is of the 
form $x^5+y^5=2a_4z^5$ and it is excluded by Lemma \ref{ell5}. Thus $6\nmid a_4$.

We consider different cases of $a_4.$  Many cases are excluded using Lemma \ref{ell5} by forming suitable  $[p,q,r]$ equations which lead to equations of the form $x^5+y^5=cz^5$ with $P(c)\leq 5.$ In some cases we use Lemma \ref{l^2}. The values of $p,q,r$ when 
Lemma \ref{ell5} is used and the values of $j_1,j_2$ when Lemma \ref{l^2} is used with $\ell=5$ and $q=5^2$ are given in Table 3. Here $j_q$ is taken as $1$ or $6$ according as $5^2|a_1$ or $a_6.$

\begin{table}
\begin{center}
\begin{tabular}{|c|c|c|}
\hline
Cases& Sub-cases& $[p,q,r]$\\
\hline
$2|a_4,3\nmid a_4$ & $\alpha_4=4$ & $[3,4,5]$\\
-& $\alpha_4=1$ & $[4,5,6]$\\
\hline
$2\nmid a_4,3|a_4$ & $5\nmid a_1a_6,\alpha_3=3,\beta_4=4$ &$[3,4,6]$\\
-&$5\nmid a_1a_6, \alpha_3=3,\beta_4=1$ &$[1,3,6]$\\
-& $5|a_1a_6,\alpha_3\in \{1,3\}$ & $(j_1,j_2)=(3,5)$\\
\hline
$2\nmid a_4,3\nmid a_4$& $\alpha_1=1$ & $[3,4,5]$\\
-&$ 5\nmid a_1a_6,\alpha_1=2$ & $[1,3,5]$\\
-&$ 5\mid a_1a_6,\alpha_1=2$ & $(j_1,j_2)=(4,5)$\\
\hline
\end{tabular}
\caption{\label{tab3:table-name} {$(k,\ell,i)=(6,5,2)$}}
\end{center}
\end{table}

\subsection{\bf Let  $(k,\ell)=(6,5)$ and $i=3$.}  We have either 
$$\al_j=0\ for\ all\  j\neq \textcolor{\defc}{3}$$ 
or 
$$
\begin{matrix}
\al_1 & \al_5 & \al_2+\al_4+\al_6 \\
4 & 1 & 0 \\ 
1 & 4 &0
\end{matrix}
\quad {\rm or} \quad 
\begin{matrix}
\al_2 & \al_4 & \al_6 & \al_1+\al_5 \\
1 & 3 & 1 &0 \\
2 & 1 & 2  &0
\end{matrix} 
$$
and either $\beta_j=0$ for all $j\neq 3$ or 
$$
\begin{matrix}
\beta_1 & \beta_4 & \sum_{j\neq 1, 4} \beta_j\\
1 & 4 & 0\\
4 & 1  &0
\end{matrix}
\quad {\rm or} \quad 
\begin{matrix}
\beta_2 & \beta_5 & \sum_{j\neq 2, 5} \beta_j\\
1 & 4 & 0 \\
4 & 1 &0.
\end{matrix} 
$$
Note that  $\beta_6=0,P(a_5)\leq 3$ and if 5 divides any $a_j,$ then $5|a_1,a_6.$   Also  by Lemma \ref{k5}, we have either 
$$3|a_1a_4\ \text{or} \ 3|a_2a_5\ \text{if}\ 5\nmid a_1a_6.$$ 

Let, $2\mid a_5$, then $(\alpha_1,\alpha_5)=(1,4),(4,1).$ We consider the equation $[1,4,5]$, i.e., $a_1x_1^5+3a_5x_5^5=4a_4x_4^5$. It is easy to conclude that, $2\mid x_1$ and $2\mid x_4$. Hence, $2\mid (n+d)$ and $2 \mid (n+4d)$. Thus $2\mid d$, and so $2\mid n$. This contradicts $\gcd(n,d)=1$. Hence $2\nmid a_5.$

In Table 4, we give various possibilities of $(a_1,a_2,a_4,a_6)$ when 2 does not divide $a_5.$ These are excluded either by considering  $[p,q,r]$ equations and using  Lemma \ref{ell5} or 
by using Lemma \ref{l^2} with $\ell=5,q=5^2$ and suitable $j_1,j_2.$ The choices  of $p,q,r$ when Lemma \ref{ell5} is used and $j_1,j_2$ when Lemma \ref{l^2} is used  
are given in the last column of the table. Here $j_q$ is taken as 1 or 6 according as $5^2|a_1$ or $a_6.$

\begin{table}
\begin{center}
\begin{tabular}{|c|c|c|}
\hline
$2\nmid a_5,3|a_5$& $5\nmid a_1a_6,\beta_2=4$ &$[1,2,5]$\\
-& $5|a_1a_6,\beta_2=4$& $(j_1,j_2)=(2,4)$\\
-& $5\nmid a_1a_6,\beta_2=1$ &$[1,4,5]$\\
-& $5|a_1a_6,\gamma_6=4,\beta_2=1$& $(j_1,j_2)=(4,5)$\\
-& $5|a_1a_6,\gamma_6=1,\beta_2=4$& $(j_1,j_2)=(2,4)$\\
\hline
$2\nmid a_5,3\nmid a_5$& $5\nmid a_1a_6$ &$[2,4,6]$\\
-& $5|a_1a_6,\gamma_6=1$& $(j_1,j_2)=(2,5)$\\
-& $5|a_1a_6,\gamma_6=4,a_2=1\ or \ 2$& $(j_1,j_2)=(2,5)$\\
-& $5|a_1a_6,\gamma_6=4,a_2=2^2,a_4=2$ &$[4,5,6]$\\
-& $5|a_1a_6,\gamma_6=4,a_2=2^2,a_4=2\cdot 3\ or\ 2\cdot 3^4$& $(j_1,j_2)=(4,5)$\\
\hline
\end{tabular}
\caption{\label{tab4:table-name} {$(k,\ell,i)=(6,5,3)$}}
\end{center}
\end{table}

\section{\bf The case $(k, \ell)=(7, 3)$}  

Let $(k, \ell)=(7, 3)$. Then $i\in\{2,3,4\}.$ Suppose 
$5\nmid a_j$ for any $j\neq i.$ By Lemma  \ref{k5} with 
$k=4, m=n+3d, s=d$ when $i\in \{2, 3\}$ and by Lemma \ref{GHP} with $m=n+d, s=d$ when 
$i=4$, we get a contradiction. Thus 
$$5|a_1a_6\ \textrm { or }\ 5|a_2a_7.$$ 

\subsection{\bf Let $(k,\ell)=(7,3)\ {\rm and}\ i=2$.} Then $5|a_1a_6$. By Lemma  \ref{k5} with $k=3, m=n+2d, s=d$, we 
may assume that 
$$3|a_3a_4a_5.$$  
Hence we have either $\al_j=0$ for all $j\neq 2$ or 
$$
\begin{matrix}
\al_1 & \al_3 & \al_5 & \al_7 & \sum_{j\neq 1, 3, 5, 7}\al_j \\
2 & 1 & 2 & 1 & 0\\
1 & 2 & 1  & 2 & 0
\end{matrix}
\quad {\rm or} \quad 
\begin{matrix}
\al_4 & \al_6 & \sum_{j\neq 4, 6}\al_j \\
2 & 1 &  0\\
1 & 2  & 0
\end{matrix} 
$$
and   
$$
\begin{matrix}
\beta_1 & \beta_4 & \beta_7  & \sum_{j\neq 1, 4, 7}\beta_j \\
1 & 1 & 1 & 0 
\end{matrix}
\quad {\rm or} \quad 
\begin{matrix}
\beta_3 & \beta_6 & \sum_{j\neq 3, 6}\beta_j\\
2 & 1 & 0 \\
1 & 2 & 0.
\end{matrix} 
$$
We have $\beta_5=0$ always.  Further $(\gamma_1, \gamma_6)\in \{(1, 2), (2, 1)\}$. 
For various possibilities of the $a_i'$s we consider 
 equations $[p, q, r]$ and exclude them by 
 Lemma \ref{ell3}. The details are given in Table 5. 

\begin{table}
\begin{tabular}{|c|c|c|c|}
\hline
Cases& Sub-Cases& $[p,q,r]$& $ab$\\
\hline
$2\nmid a_j$& $3|a_1a_4a_7$&$[4,5,7]$ &$2$\\
-& $\beta_3=2$&$[3,5,7]$ &18\\
-& $\beta_3=1,\gamma_1=1$&$[1,4,7]$ &$10$\\
-& $\beta_3=1,\gamma_1=2$&$[4,6,7]$ &$10$\\
\hline
$2|a_j$& -&- &-\\
$2|a_1a_3a_5a_7$&$3|a_1a_4a_7$ & $[4,5,7]$& $2$\\
-& $3|a_3a_6$, $(a_3,a_7)\neq(12,4)$ & $[3,4,7]$ & $\{1,2,18\}$\\
-& $3|a_3a_6$, $(a_3,a_7)=(12,4)$ & $[3,5,7]$ & $3$\\
$2|a_4a_6$& $3|a_1a_4a_7$ & $[4,5,7]$ & $\{1,4\}$\\
-& $3|a_3a_6$, $(a_3,a_4)\neq(2,3)$ & $[3,4,7]$ & $\{1,2,18\}$\\
-& $3|a_3a_6$, $(a_3,a_4)=(2,3)$ & $[1,6,7]$ & $5$\\
\hline
\end{tabular}
\caption{\label{tab5:table-name} {$(k,\ell,i)=(7,3,2)$}}
\end{table}

\subsection{\bf Let $(k,\ell)=(7,3)\ {\rm and}\ i=3$} 
By Lemma \ref{k5} with $k=4, m=n+3d,s=d$ we may assume that
$$5|a_1a_6\ \textrm {or }\ 5|a_2a_7.$$
And $(\gamma_1, \gamma_6)\in \{(1, 2), (2, 1)\}$ or 
$(\gamma_2, \gamma_7)\in \{(1, 2), (2, 1)\}.$
Suppose 2 and 3 do not divide any $a_j.$ Then
$$a_2=a_4=a_5=a_7=1\ \textrm { if }\ 5|a_1a_6$$
or
$$ a_1=a_4=a_5=a_6=1\ \textrm { if }\ 5|a_2a_7.$$
The first possibility is ruled out by Lemma \ref{ell3} using
$$[5,6,7] \ \textrm { with } \ ab=10\ \textrm { if } \ 5^2||a_1\ \textrm { and }\ 5||a_6$$
and
$$[1,2,5]\ \textrm { with }\ ab=60\ \textrm { if }\ 5||a_1\ \textrm { and } \ 5^2||a_6.$$
The second possibility is excluded by $[4,5,6]$ with $ab=2.$
Hence we may assume that either 
$$2|a_j \ \textrm { or }\ 3|a_j.$$
Further we have either
 $$
\begin{matrix}
\al_1 & \al_5 & \al_7 & \sum_{j\neq 1, 5, 7}\al_j \\
0 & 2 & 1 & 0\\
2 & 0 & 1  & 0 \\
1 & 1 & 1  & 0 \\
\end{matrix}
\quad {\rm or} \quad 
\begin{matrix}
\al_2 & \al_4 & \al_6 & \sum_{j\neq 2, 4, 6}\al_j \\
0 & 1 &  2\\
1 & 1  & 1\\
2 & 1 & 0
\end{matrix} 
$$
whenever $2|a_j$ and   
$$
\begin{matrix}
\beta_1 & \beta_4 & \beta_7  & \sum_{j\neq 1, 4, 7}\beta_j \\
1 & 1 & 1 & 0 
\end{matrix}
\quad {\rm or} \quad 
\begin{matrix}
\beta_2 & \beta_5 & \sum_{j\neq 2, 5}\beta_j\\
2 & 1 & 0 \\
1 & 2 & 0
\end{matrix} 
$$
when $3|a_j.$
Let $5|a_2a_7$. By Lemma \ref{k5} with $k=3, m=n+3d, s=d$, we may assume that 
$3|a_4a_5a_6$ and hence 
$$3|a_4a_5$$ 
since $3\nmid a_6$. We will use this while listing the case $5|a_2a_7.$
In Tables 6 and 7, we give the various possibilities according as 
$5|a_1a_6$ or $5|a_2a_7$ respectively and exclude them by considering equations $[p, q, r]$ along with Lemma \ref{ell3}. 

\begin{table}
\begin{center}
\begin{tabular}{|c|c|c|c|}
\hline
Cases & Sub-Cases & $[p,q,r]$ & $ab$\\
\hline
$2\nmid a_j$ & $\gamma_1=1$ & $[1,4,7]$ & $10$\\
-& $\gamma_1=2, 3|a_1a_4a_7$ & $[4, 6, 7]$  &  $10$ \\ 
-& $\gamma_1=2, 3|a_2a_5$ & $[1, 2, 7]$ &  $ \{2, 18\}$ \\ 
 \hline
 $3\nmid a_j$  & $\gamma_1=1$ & $[2, 6, 7]$ &  $\{1, 2, 4\}$  \\
 -& $\gamma_1=2$ & $[1, 5, 7]$ &  $150$  \\ 
 \hline
$3|a_1a_4a_7$ & $\gamma_1=1$ & $[1, 5, 7]$  &  $10$    \\ 
-& $\gamma_1=2, j_2\neq 5$ & $[1, 4, 7]$ &  $\{25, 100\}$    \\ 
-& $\gamma_1=2, j_2=5$ & $[2,4,5]$  &  $18$  \\ 
\hline

$3|a_2a_5$ & $\gamma_1=1, \beta_2=2$ & $[2, 6, 7]$  &  $\{18, 36\}$    \\ 
-& $\gamma_1=1, \beta_2=1$ & $[1, 5, 7]$  &  $10$    \\ 
-& $\gamma_1=2, \beta_2=1, (a_1,a_2,a_7)\neq (50,3,2)$ & $[1, 2, 7]$  &  $ \{18, 36\}$    \\
-& $\gamma_1=2, \beta_2=1, (a_1,a_2,a_7)= (50,3,2)$ & $[4, 6, 7]$  &  $60$    \\ 
-& $\gamma_1=2, \beta_2=2$ & $[1, 2, 7]$  &  $\{1, 2, 4\}$    \\ \hline
\end{tabular}
\end{center}
\caption{\label{tab6:table-name} {$(k,\ell,i)=(7,3,3),5|a_1a_6$}}
\end{table}
\begin{table}
\begin{center}
\begin{tabular}{|c|c|c|c|} \hline
Cases &Sub-cases&  $[p,q,r]$ & $ab$\\
\hline
$3|a_2a_5$ & $\beta_2=1, (a_4,a_5,a_6)\neq (1,36,1),(2,9,2)$&$[4,5,6]$& $\{18,36\}$\\
- & $\beta_2=1, (a_4,a_5,a_6)= (1,36,1)$&$[4,6,7]$& $\{60, 300\}$\\
- & $\beta_2=1, (a_4,a_5,a_6)= (2,9,2)$&$[4,5,7],[2,4,6]$ & $100,150$\\
- &$\beta_2=2, 2\nmid a_2a_4a_6$,  &$[1,4, 5]$&  $\{18, 36\}$ \\
- &$\beta_2=2, 2\mid a_2a_4a_6$ &$[2,4, 5], [5,6,7]$&  $\{60, 150, 300\}$ \\
\hline
$3|a_1a_4a_7$ &$\gamma_2=1$& $[1, 6,7]$&  $ \{1,2,4\}$ \\ 
-& $\gamma_2=2$&$[2,4,6]$&$150$\\
\hline
\end{tabular}
\end{center}
\caption{\label{tab7:table-name} {$(k,\ell,i)=(7,3,3),5|a_2a_7$}}
\end{table}
\subsection{\bf Let $(k,\ell)=(7,3)\ {\rm and}\ i=4$.} Then $5|a_1a_6$ or $5|a_2a_7$. By Lemma \ref{k5} 
with $k=3, m=n+d, s=2d$ if $5|a_1a_6$ and $k=3, m=n-d, s=2d$ if $5|a_2a_7$, 
we may assume that 
$$3|a_1a_2a_3.$$  
Also we have either $\al_j=0$ for all $j\neq 2$ 
or 
$$
\begin{matrix}
\al_1 & \al_3 & \al_5 & \al_7 & \sum_{j\neq 1, 3, 5, 7}\al_j \\
2 & 1 & 2 & 1 & 0\\
1 & 2 & 1  & 2 & 0
\end{matrix}
\quad {\rm or} \quad 
\begin{matrix}
\al_2 & \al_6 & \sum_{j\neq 2, 6}\al_j \\
2 & 1 &  0\\
1 & 2  & 0
\end{matrix} 
$$
and   
$$
\begin{matrix}
\beta_1 & \beta_7  & \sum_{j\neq 1, 7}\beta_j \\
2 & 1 & 0 \\
1 & 2 & 0 
\end{matrix}
\quad {\rm or} \quad 
\begin{matrix}
\beta_2 & \beta_5 & \sum_{j\neq 2, 5}\beta_j\\
2 & 1 & 0 \\
1 & 2 & 0
\end{matrix} 
\quad {\rm or} \quad 
\begin{matrix}
\beta_3 & \beta_6 & \sum_{j\neq 3, 6}\beta_j\\
2 & 1 & 0 \\
1 & 2 & 0.
\end{matrix} 
$$
Further $(\gamma_1, \gamma_6)\in \{(1, 2), (2, 1)\}$ or 
$(\gamma_2, \gamma_7)\in \{(1, 2), (2, 1)\}$. In Tables 8 and 9, we list the possibilities when $5|a_1a_6$ and $5|a_2a_7$ and exclude these by 
 Lemma \ref{ell3} by considering equations
$[p, q, r]$ which are of the form $x^3+ay^3+bz^3=0$ except the case $(a_1,a_2,a_3,a_5,a_6,a_7)=(10,3,4,18,25,4)$. In this case, $n=-17$ and $d=7$, gives rise to a solution.\\
\begin{table}
\begin{center}
\begin{tabular}{|c|c|c|c|} \hline
Cases & Sub-cases & $[p,q,r]$ & $ab$ \\
\hline
$9|a_3a_5a_7$ & $9|a_3$& $[3, 5, 7], [1, 3, 7]$  &  $\{18, 36\}, \{5, 25\}$ \\
- & $9|a_5$& $[3, 5, 7], [1, 2, 7]$  &  $\{18, 36\}, 18$ \\
- & $9|a_7$& $[3, 5, 7], [1, 2, 7], [1, 2, 3]$  &  $\{18, 36\}, 150$ \\
\hline
$9|a_6$ & - & $[2, 3, 5], [1, 3, 6]$  &  $ \{18, 36\}, 60$ \\  \hline

$9|a_1$ & $ \al_6\neq 1$ & $[1, 3, 6]$,  &  $\{5, 10\}$ \\ 
 & $ \al_6=1$ & $[1, 5, 6]$  &  $45$ \\   \hline

$9|a_2$ & $\gamma_1=2$ & $[1, 2, 7]$,  &  $ \{1, 2, 4\}$ \\  
 & $ \gamma_1=1$ & $[2, 6, 7], [1, 3, 7]$,  &  $ \{18, 36\}, 60$ \\   \hline
\end{tabular}
\end{center}
\caption{\label{tab8:table-name} {$(k,\ell,i)=(7,3,4), 5|a_1a_6$}}
\end{table}

\begin{table}
\begin{center}
\begin{tabular}{|c|c|c|c|} \hline
Cases & Sub-Cases & $[p,q,r]$ & $ab$\\
\hline
$9|a_1a_3a_5$ &- &$[1, 3, 5], [1, 2, 5],$  &  $\{18, 36\}, 25, $ \\
- &- &$[1, 2, 6], [1, 2, 3], [1, 2, 7], [5, 6, 7]$  &  $150$ \\
\hline
$9|a_2$ &  & $[3, 5, 6], [3, 6, 7], [1, 2, 6]$ &  $\{18, 36\}, 150, 36$ \\  \hline
$9|a_6$ & - & $[1, 6, 7]$  &  $\{1, 2, 4, 25, 100\}$ \\ 
&  - & $[1, 2, 3]$  &  $150$ \\  \hline

$9|a_7$ & $ \al_2\neq 1$ & $[2, 5, 7]$  &  $\{5, 10\}$ \\ 
& $ \al_2=1$ & $[1, 2, 5]$  &  $\{45, 225\}$ \\  \hline
\end{tabular}
\end{center}
\caption{\label{tab9:table-name} {$(k,\ell,i)=(7,3,4),5|a_2a_7$}}
\end{table}

\section{\bf The case $(k, \ell)=(8, 7)$}  
We have  $i\in\{2,3,4\}.$
\subsection{\bf Let  $(k,\ell)=(8, 7)$ and $i=2$.}  
We have either $7|a_1a_8$ or $7\nmid a_j$ for any $j\neq 2$. By 
Lemma \ref{k5} with $k=5, m=n+2d, s=d$, we may assume that $5|a_j$ 
for some $j$. Then either 
$$5|a_1a_6\ \textrm{ or } \ 5|a_3a_8.$$ 
Also $3|a_j$ for some $j$ by 
considering Lemma \ref{k5} with $k=3$ and $(m, s)=(n+2d, d), (n+3d, d)$ according as 
$5|a_1a_6$ or $5|a_3a_8$, respectively. Hence we have either 
$\al_j=0\ \textrm{ for all }\ j\neq 2$ 
or 
$$
\begin{matrix}
\al_1 & \al_3 & \al_5 & \al_7 & \sum_{j\neq 1, 3, 5, 7}\al_j \\
3 & 1 & 2 & 1 & 0\\
1 & 3 & 1  & 2 & 0\\
2 & 1 & 3  & 1 & 0\\
1 & 2 & 1  & 3 & 0\\
\end{matrix}
\quad {\rm or} \quad 
\begin{matrix}
\al_4 & \al_6 & \al_8 & \sum_{j\neq 4, 6, 8}\al_j \\
4 & 1 & 2 & 0\\
1 & 5 & 1  & 0\\
2 & 1 & 4  & 0
\end{matrix} 
$$
and   
$$
\begin{matrix}
\beta_1 & \beta_4 & \beta_7  & \sum_{j\neq 1, 4, 7}\beta_j \\
5 & 1 & 1 & 0 \\
1 & 5 & 1 & 0 \\
1 & 
1 & 5 & 0 
\end{matrix}
\quad {\rm or} \quad 
\begin{matrix}
\beta_3 & \beta_6 & \sum_{j\neq 3, 6}\beta_j\\
6 & 1 & 0 \\
1 & 6 & 0.
\end{matrix} 
\quad {\rm or} \quad 
\begin{matrix}
\beta_5 & \beta_8 & \sum_{j\neq 5, 8}\beta_j\\
6 & 1 & 0 \\
1 & 6 & 0.
\end{matrix} 
$$

Suppose $7\nmid a_1a_8$.  The case $5|a_1a_6$ is excluded by Lemma \ref{GHP} with 
$m=n+3d$ and $s=d$. Hence 
$$5|a_3a_8.$$
 If $3\nmid a_1a_4a_7$, then Lemma \ref{k5} 
with $k=3, m=n-2d, s=3d$ gives a contradiction. Thus 
$$3|a_1a_4a_7$$ 
so that $3\nmid a_j$ 
for $j\in \{3, 5, 6, 8\}$. We consider $[5, 6, 3, 8]$ which gives  
\begin{align*}
(a_5a_6)(x_5x_6)^7-(a_3a_8)(x_3x_8)^7=6d^2. 
\end{align*}
Note that $2\nmid a_j$ for any $j$ if $2|d$.  
We get an equation of the form Lemma \ref{eq4} $(i), (iii)$, or Lemma \ref{ell7} $(iii)$ according as $2\nmid d$ or 
$2|d$, respectively which is excluded by Lemma \ref{eq4}-\ref{ell7}.  
Thus 
$$7|a_1a_8.$$ 
Then $j_{7^2}\in \{1, 8\}$ and by \eqref{chi} in Lemma \ref{l^2}, we have 
$$\left(\frac{a_4a_5}{12}\right)^6\equiv \left(\frac{a_3a_6}{10}\right)^6 \pmod{7^2}.$$
Since $5|a_1a_6$ or $5|a_3a_8$, we have $\nu_5(a_3a_6)=1+\gamma$ with $\gamma\in \{0, 5\}$. 
Also note that $3\nmid a_3a_6$, since $3\mid a_1a_4a_7$.
Hence, we write $a_3a_6=2^{\al_3+\al_6} 5^{1+\gamma}$. 
Hence, using the above equivalence and observing that 
$\La_{7^2}=\{\pm  1,\pm   18,\pm   19\}$, we obtain 
\begin{align*}
\frac{a_4a_5}{12}\equiv \frac{\la 2^{\al_3+\al_6}5^{\gamma}}{2} \pmod{7^2} \quad {\rm with}  
\quad \gamma\in \{0, 5\} \ {\rm and} \ \la\in \La_{7^2}. 
\end{align*}
Further writing $a_4a_5=2^{\al_4+\al_5}3^{\beta_4+\beta_5}$ and using 
$$\frac{1}{6}\equiv -8, \quad \la 5^5\equiv \pm  2,\pm   11,\pm   13 \pmod{7^2}, $$
we obtain
\begin{align}\label{45-36}
2^{3+\al_4+\al_5-\al_3-\al_6}3^{\beta_4+\beta_5}\equiv \chi  \pmod{7^2}
\end{align}
where 
\begin{align*}
\chi \in \{ \pm 1, \pm  18, \pm  19\} \quad & {\rm if} \ \nu_5(a_3a_6)=1,\\ 
\chi \in \{ \pm 2, \pm  11, \pm  13\} \quad & {\rm if} \ \nu_5(a_3a_6)=6.
\end{align*}
Observe that $\beta_4+\beta_5\in \{0, 1, 5, 6\}$. For each $\beta\in \{0, 1, 5, 6\}$ we check for 
the possibilities of $3+\al_4+\al_5-\al_3-\al_6$ satisfying \eqref{45-36}. We use 
$2^{-1}\equiv 25 \pmod{49}$. These are listed in Table 10.
\begin{table}
\begin{center}
\begin{tabular}{|c|l|l|}
\hline
$\beta_4+\beta_5=0$ & $\al_3=3, \al_5=1, \gamma=5$; & $a_4=1, a_5=2$  \\  \hline 
$\beta_4+\beta_5=1$ & $\al_3=1, \al_5=2, \gamma=0$; & $a_4=1, a_5=12$  or  $a_4=3, a_5=4$  \\ \hline
- & $\al_4=2, \al_6=1, \gamma=0$;  & $a_4=4, a_5=3$  \ or $a_4=12, a_5=1$  \\  \hline
 - &  $\al_3=1, \al_5=3, \gamma=5$;  & $a_4=1, a_5=24$  or $a_4=3, a_5=8$  \\ \hline 
$\beta_4+\beta_5=5$ & $\al_4=1, \al_6=5, \gamma=0$;  & $a_4=2\cdot 3^5, a_5=1$ \\  \hline   
 $\beta_4+\beta_5=6$ & $\al_3=1, \al_5=2, \gamma=5$;  & $a_4=1, a_5=4\cdot 3^6$  \\  \hline
 -& $\al_4=2, \al_6=1, \gamma=5$; & $a_4=4, a_5=3^6$  \\ \hline 
\end{tabular}
\caption{\label{tab10:table-name} {$(k,\ell,i)=(8,7,2)$}}
\end{center}
\end{table}
Taking $j_{7^2}\in \{1, 8\}$ and using $(a_4/(4-j_{7^2}))^6\equiv (a_5/(5-j_{7^2}))^6 \pmod{7^2}$ 
by \eqref{chi} in Lemma \ref{l^2} and further using 
$$\left(\frac{a_4}{4-j_{7^2}}\right)^6\equiv \left(\frac{a_j}{j-j_{7^2}}\right)^6 \pmod{7^2} \ {\rm with} \  
\begin{cases}
j=7 & {\rm when} \ \beta_4+\beta_5=0\\
j=3 & {\rm when} \ \beta_4+\beta_5\neq 0,
\end{cases}
$$
the cases in Table 10 are excluded except those listed in Table 11. 
\begin{table}
\begin{center}
\begin{tabular}{|c|c|c|l|}
\hline
1) & $\beta_4+\beta_5=1$ & $j_{7^2}=1, \gamma=0$; & $ a_1=2^33^55^67^6, a_j=j-1$ for $3\leq j\leq 8$  \\  \hline

2) &- & $j_{7^2}=8, \gamma=0$; & $a_8=2^33^65^67^6, a_j=8-j$ for $j\neq 2, 8$  \\  \hline
3)&-& $j_{7^2}=1, \gamma=5$;  & $ a_3=2\cdot 5^6, a_5=24, a_6=1, a_8=3\cdot 5\cdot 7$   \\ \hline 
 4)& $\beta_4+\beta_5=6$ & $j_{7^2}=1, \gamma=5$; & $a_1=2^3\cdot 7^6, a_3=2\cdot 5^6, a_5=4\cdot 3^6$,\\
 &&& $ a_7=2, a_4=a_6=1, a_8=3\cdot 5\cdot 7.$  \\  \hline 
\end{tabular}
\caption{\label{tab11:table-name} {$(k,\ell,i)=(8,7,2)$}}
\end{center}
\end{table}
By considering $[5, 6, 3, 8]$, the possibilities $2)$ and $3)$ are excluded. For  
$4)$, we use Lemma \ref{l^2} with $\ell=7$ and $q=29$.  
Note that $\la_{29}\in \{ \pm 1, \pm  12\}.$ Suppose $29|d.$ Then by taking $(j_1,j_2)=(4,7)$ in \eqref{q|d}, we get a contradiction. Let $29\nmid d.$ Suppose $j_{29}\neq 4,7.$ Using \eqref{chi} with $(j_1,j_2)=(4,7),$ we see that $j_{29}\in \{1,5,19,25\}.$ These cases are excluded using \eqref{chi} with $(j_1,j_2)=(5,7).$
When $j_{29}=4\ \textrm {or }\ 7$, then we use \eqref{chi} with $(j_1,j_2)=(6,7)$ or $(4,6)$ to get a contradiction. Finally we consider 
$1).$ Using $[3,4,7]$
we see that
$$x_3^7+x_7^7=2x_4^7$$
which has no solution  by Lemma \ref{eq1} with $p=2.$

\subsection{\bf Let  $(k,\ell)=(8, 7)$ and $i=3$.}  
We have either $7|a_1a_8$ or $7\nmid a_j$ for any $j\neq 3$.
Also either 
$5|a_1a_6$ or $5|a_2a_7$ or $5\nmid a_j$ for any $j\neq 3$. Further 
we have either 
$\al_j=0\ \textrm { for all } \ j\neq 3$
or 
$$
\begin{matrix}
\al_1 & \al_5 & \al_7 & \sum_{j\neq 1, 5, 7}\al_j \\
4 & 2 & 1 & 0\\
2 & 4 & 1  & 0\\
1 & 1 & 5  & 0
\end{matrix}
\quad {\rm or} \quad 
\begin{matrix}
\al_2 & \al_4 & \al_6 & \al_8 & \sum_{j\neq 2, 4, 6, 8}\al_j \\
3 & 1 & 2 & 1 & 0\\
1 & 3 & 1  & 2 & 0\\
2 & 1 & 3  & 1 & 0\\
1 & 2 & 1  & 3 & 0\\
\end{matrix} 
$$
and   
$$
\begin{matrix}
\beta_1 & \beta_4 & \beta_7  & \sum_{j\neq 1, 4, 7}\beta_j \\
5 & 1 & 1 & 0 \\
1 & 5 & 1 & 0 \\
1 & 1 & 5 & 0 
\end{matrix}
\quad {\rm or}  \quad 
\begin{matrix}
\beta_2 & \beta_5 & \beta_8 & \sum_{j\neq 2, 5, 8}\beta_j\\
5 & 1 & 1 & 0 \\
1 & 5 & 1 & 0 \\
1 & 1 & 5 & 0.
\end{matrix} 
$$
Note that $\beta_6=0$. Let $5|a_2a_7$. By considering Lemma \ref{k5} with $k=3$ and 
$(m, s)=(n+3d, d)$, we may assume that 
$$3|a_4a_5.$$  
We have  
$$5^7|a_2a_7\ \textrm { and } \ 7^\delta|a_1a_8,\delta\in \{0, 7\}$$
with  
$$P(a_2a_7/5^7)\leq 3\ \textrm{ and }\ P(a_1a_8/7^\delta)\leq 3.$$
 We consider different 
equations as follows:
\begin{align*}
[4, 5, 1, 8] &\quad  {\rm if} \quad \beta_1=\beta_4=1 \ {\rm or} \ \beta_5=\beta_8=1\\
[2, 7, 1, 8] &\quad  {\rm if} \quad \beta_1=\beta_7=1 \ {\rm or} \  \beta_2=\beta_8=1\\
[4, 5, 2, 7] &\quad  {\rm if} \quad \beta_4=\beta_7=1 \ {\rm or} \ \beta_2=\beta_5=1
\end{align*}
When $2|d$, we have all $a_j'$s odd and above resulting equations are of the form 
$x^7+y^7=cz^2$ with $c\in \{1, 2\}$  and hence excluded by Lemma \ref{eq4} $(iii)$. 
When $d$ is odd, we take \emph{   }[4, 5, 2, 7] which has the form as in Lemma \ref{eq4} $(iv)$ 
and hence excluded. Thus 
$$5\nmid a_2a_7.$$ 

Suppose $7\nmid a_1a_8$. By Lemma \ref{k5} with $k=5, m=n+3d, s=d$, we may 
assume that $5|a_j$ for some $j$. Then $5|a_1a_6$. This is excluded by Lemma \ref{GHP} 
with $m=n+3d$ and $s=d$. 
Thus 
$$7|a_1a_8.$$ 
Then $j_7\in \{1, 8\}$ and by \eqref{chi} in Lemma \ref{l^2} with $q=7^2$, we have 
$$\left(\frac{a_4a_5}{12}\right)^6\equiv \left(\frac{a_2a_7}{6}\right)^6 \pmod{7^2}.$$
Hence we obtain
\begin{align}\label{45-27}
2^{\al}3^{\beta}\equiv \la  \pmod{7^2} \quad {\rm with }\ \la \in \{\pm  1, \pm  18, \pm  19\}
\end{align}
where 
\begin{align*}
\al=1+\al_2-\al_4+\al_7-\al_5 \quad {\rm and} \quad \beta=\beta_2-\beta_5+\beta_7-\beta_4.
\end{align*}
Observe that $\beta\in \{0,   4\}$ and $-2\leq \al\leq 5$. Hence \eqref{45-27} implies 
$$(\al, \beta)\in \{(0, 0), (2, 4), (5, -4), (-2, -4)\}.$$
These lead to the possibilities as shown in the first two columns in Table 12. 
\begin{table}
\begin{center}
\begin{tabular}{|c|l|l|}
\hline
$\al$ & $\beta$  & $a_{j_1}, a_{j_2}$ \\  \hline 
$\al=0;$ & $\beta=\beta_2=\beta_4=0$ &  $a_5=a_7=1$   \\  
$\al_2=1, \al_4=2$& $\beta=0, \beta_4=\beta_7=1$   & $a_5=1, a_7=3$   \\ 
&$\beta=0, \beta_2=\beta_5=1$ & $a_5=3, a_7=1$   \\ \hline

$\al=0;$ & $\beta=\beta_2=\beta_4=0$ & $a_2=a_4=1$  \\ 
$\al_5=2, \al_7=1$ & $\beta=0, \beta_4=\beta_7=1$ &  $a_2=1, a_4=3$ \\ 
 & $\beta=0, \beta_2=\beta_5=1$ & $a_2=3, a_4=1$  \\ \hline 

$\al=2;$ & $\beta=4, \beta_4=1, \beta_7=5$ & $a_2=4, a_5=1$   \\ 
$\al_2=2, \al_4=1$ & $\beta=4, \beta_2=5, \beta_5=1$ & $a_4=2, a_7=1$   \\ \hline

$\al=5;$ & $\beta=-4, (\beta_4, \beta_7)=(5, 1)$ & $a_2=1, a_5=2$  \\ 
$\al_5=1, \al_7=5$ & $\beta=-4, (\beta_2, \beta_5)=(1, 5)$ & $a_2=3, a_4=1$  \\ \hline

$\al=-2;$ & $\beta=-4, (\beta_4, \beta_7)=(5, 1)$ & $a_2=1, a_5=2^4$  \\ 
$\al_5=4, \al_7=1$ & $\beta=-4, (\beta_2, \beta_5)=(1, 5)$ & $a_2=3, a_4=1$  \\ \hline
\end{tabular}
\end{center}
\caption{\label{tab12:table-name} {$(k,\ell,i)=(8,7,3)$}}
\end{table}
All the possibilities are excluded by \eqref{chi} of Lemma \ref{l^2} with $(a_{j_1},a_{j_2})$ as indicated in the third column with $ q=7^2$ and $j_q\in \{1,8\}$ except the following cases.
 
\begin{align*}
(1)\ &a_1=7, a_2=6, a_4=4, a_5=3, a_6=2, a_7=1, a_8=2^3\cdot 3^5\cdot 7^6;\\
(2) \ &a_2=6, a_4=4, a_5=3, a_7=1, a_8=2^3\cdot 3^5\cdot 7^6, \\
\ &(a_1,a_6)\in \{(7\cdot 5,2\cdot 5^6),( 7\cdot 5^6,10)\};\\
(3) \ &a_2=1, a_4=3, a_5=4, a_7=6, a_8=7, a_6=1, a_1=2^4\cdot 3^5\cdot 7^6;\\
(4) \ &a_2=1, a_4=3, a_5=4, a_7=6, a_8=7, \\
           \ &   (a_1,a_6)\in \{(2^4\cdot 3^5\cdot 5^6\cdot 7^6,5),(2^4\cdot 3^5\cdot 5\cdot 7^6, 5^6)\}.
\end{align*}
The possibility $(2)$ and $(3)$ are excluded by  taking $\ell=7, q=7^2, (j_1, j_2)=(2, 6)$ in
\eqref{chi} with $j_{7^2}=8, 1$, respectively.  For $(1),$ consider 
$[2,5,6]$ to get
$x^7_6+x^7_2=2x^7_5$ which has no solution by Lemma \ref{eq1} with $p=2.$ 
 For $(4),$ consider  $[2, 5, 1, 6]$  which is of the form $x^7+2^2\cdot 3^5\cdot 7^6y^7=z^2$ 
with $5|xy$. This is excluded by Lemma \ref{ell7} $(ii)$. 
\subsection{\bf Let  $(k,\ell)=(8, 7)$ and $i=4$.} We have either $5|a_ja_{5+j}$ for some 
$j\in \{1, 2, 3\}$ or $5\nmid a_j$ for any $j\neq 4$. Also $7\nmid a_j$ for 
$2\leq j\leq 7, j\neq 4$. By Lemma \ref{GHP} with $m=n+d, s=d$, we may assume that 
 $$5|a_1a_6\ \textrm { or } \ 5|a_2a_7\ \textrm {  or }\ 5|a_3a_8.$$   
We also have either 
$$\al_j\beta_j=0\ \textrm { for } \ j\neq 4$$ 
or 
$$
\begin{matrix}
\al_1 & \al_3 & \al_5 & \al_7 & \sum_{j\neq 1, 3, 5, 7}\al_j \\
3 & 1 & 2 & 1 & 0\\
1 & 3 & 1  & 2 & 0\\
2 & 1 & 3  & 1 & 0\\
1 & 2 & 1  & 3 & 0\\
\end{matrix}
\quad {\rm or} \quad 
\begin{matrix}
\al_2 & \al_6 & \al_8 & \sum_{j\neq 4, 6, 8}\al_j \\
4 & 2 & 1 & 0\\
2 & 4 & 1  & 0\\
1 & 1 & 5  & 0
\end{matrix} 
$$
and   
$$
\begin{matrix}
\beta_1 & \beta_7 & \displaystyle{\sum_{j\neq 1, 7}}\beta_j\\
6 & 1 & 0 \\
1 & 6 & 0.
\end{matrix} 
\quad {\rm or} \quad 
\begin{matrix}
\beta_2 & \beta_5 & \beta_8  & \displaystyle{\sum_{j\neq 2, 5, 8}}\beta_j \\
5 & 1 & 1 & 0 \\
1 & 5 & 1 & 0 \\
1 & 1 & 5 & 0 
\end{matrix}
\quad {\rm or} \quad 
\begin{matrix}
\beta_3 & \beta_6 & \displaystyle{\sum_{j\neq 3, 6}}\beta_j\\
6 & 1 & 0 \\
1 & 6 & 0.
\end{matrix} 
$$
Let $5|a_2a_7$. Take the equations $[3, 6, 2, 7]$ or  $[2, 7, 1, 8]$  according as $2|d$ or $2\nmid d$, respectively. These are excluded by 
Lemma \ref{eq4} $(ii)$, $(iv)$, and Lemma \ref{ell7} $(iii).$ Thus 
$$5|a_1a_6\ \textrm { or }\ 5|a_3a_8.$$
Suppose $7\nmid a_1a_8$.
Let  $5|a_3a_8$. By Lemma \ref{k5} with 
$k=3, m=n+4d, s=d$, we may assume that 
$$3|a_5a_6a_7.$$ 
Consider equation $[5, 6, 3, 8].$ By Lemma \ref{eq4} $(iv)$ we may suppose 
that $2\nmid a_j$ for $j\neq 4$ and  $3|a_3a_6$ or $3|a_5a_8.$ This leads to an equation as in Lemma \ref{ell7} $(iii)$ which has no solution.
When $5\mid a_1a_6,$ consider equation $[2,5,1,6]$ to get a contradiction from Lemma \ref{eq4} $(iv)$ and Lemma \ref{ell7} $(iii).$
Thus 
$$7|a_1a_8.$$ 
Then $j_{7^2}\in \{1, 8\}$ and by \eqref{chi} in Lemma \ref{l^2} with 
$\ell=7, q=7^2$, we have 
\begin{align*}
 \left(\frac{a_3a_6}{10}\right)^6\equiv \left(\frac{a_2a_5}{4}\right)^6 \pmod{7^2} \quad & {\rm if} \ j_{7^2}=1;\\
 \left(\frac{a_3a_6}{10}\right)^6\equiv \left(\frac{a_2a_5}{18}\right)^6 \pmod{7^2} \quad & {\rm if} \ j_{7^2}=8.
\end{align*}
Since $5|a_1a_6$ or $5|a_3a_8$, we have $\nu_5(a_3a_6)=\gamma_3+\gamma_6\in \{1, 6\}$.  
Also note that either $3\nmid a_3a_6$ or $\nu_3(a_3a_6)=7$, and 
either $3\nmid a_2a_5$ or $\nu_3(a_2a_5)=7-\beta_8$. Hence 
$$\left(\frac{a_3a_6}{a_2a_5}\right)^6\equiv (2^{\al_6-\al_2+\al_3-\al_5} 3^{\beta_8}
5^{\gamma_3+\gamma_6})^6 \pmod{7^2}.$$ 
We observe that the solutions of 
$\la^6\equiv (10/4)^6, (10/18)^6 \pmod{7^2}$ are given by $\la\in \La_1, \la\in \La_8,$ 
respectively, where 
\begin{align*}
\La_1=\{\pm  4,\pm   22,\pm   23\} \quad {\rm and} \quad 
\La_8=\{\pm  6,\pm   10, \pm  16\}.
\end{align*}
Therefore from the above equivalences, we obtain 
\begin{align}\label{25-36}
2^{\al_6-\al_2+\al_3-\al_5} 3^{-\beta_8}
5^{\gamma_3+\gamma_6}\equiv \la  \pmod{7^2} \quad {\rm for \ some } \quad \la\in \La_1\ \textrm { or }\ \La_8. 
\end{align}
For different choices of $\beta_8\in \{0, 1, 5\}$ and $\gamma=\gamma_3+\gamma_6\in \{1, 6\}$, we have 
the possibilities as shown in Table 13.
\begin{table}
\begin{center}
\begin{tabular}{|c|c|c|c|c|}
\hline
-& $\beta_8$ & $\gamma$ & $j_{7^2}$  & -   \\  \hline 
1) & $0$ & $1$ & 1 & $\al_3=1, \al_5=2, \al_7=1$   \\  
2) & 0 & 1 & 8    & $\al_3=2, \al_5=1, \al_7=3$   \\ \hline

3) &0 & 6 & 1 &  $\al_3=1, \al_5=3, \al_7=1$  \\  
4) &0 &6& 1 &   $\al_2=2, \al_6=4, \al_8=1$   \\  
5) &0 & 6 & 8& $\al_2=\al_6=1, \al_8=5$  \\ \hline

6) & 1 & 1 & 8 &   $\al_3=1, \al_5=3, \al_7=1$   \\  
7) &1&1 & 8  & $\al_2=2, \al_6=4, \al_8=1$    \\ \hline

8)&1 &6 & 1 &  $\al_3=3, \al_5=1, \al_7=2$   \\  
9) &1&6& 1& $\al_2=4, \al_6=2, \al_8=1$  \\ \hline

10) &5 &1& 1&  $\al_3=1, \al_5=3, \al_7=1$   \\  
11) & 5& 1 & 1 &  $\al_2=2, \al_6=4, \al_8=1$  \\
12)& 5 & 1& 8& $\al_2=\al_6=1, \al_8=5$   \\ \hline

13) & 5 &6 & 8 &  $\al_3=1, \al_5=2, \al_7=1$  \\  
14) & 5&6 &8  & $\al_3=2, \al_5=1, \al_7=3$  \\ \hline 
\end{tabular}
\end{center}
\caption{\label{tab13:table-name} {$(k,\ell,i)=(8,7,4)$}}
\end{table}

Now take $\ell=7$, $q=7^2.$
In Table 14, we list the cases from Table 13 and the corresponding $(j_1,j_2)$ that are excluded by checking that \eqref{chi} is not satisfied. 
\begin{table}
 \begin{center}
\begin{tabular}{|c|c|c|c|}
\hline
 $\beta_8$ & - & Cases & $(j_1, j_2)$ \\  \hline 
0& $3|a_3a_6$ &$1)-5)$ &  $(5, 7)$   \\  \hline

$\{1, 5\}$ & $ 3|a_2a_5a_8$ & $6)-14)$   & $5|a_1a_6: (3, 7)$   \\ \hline
$\{1, 5\}$ & $3|a_2a_5a_8$ & $6)-9),11),13),14) $  & $5|a_3a_8: (6, 7)$   \\ \hline
0& $ 3|a_1a_7$ &$2)-4)$ &  $(2, 5)$   \\  \hline
\end{tabular}
\end{center}
\caption{\label{tab14:table-name} {Excluded cases for $(k,\ell,i)=(8,7,4)$}}
\end{table}
We are left with the following possibilities. 
$$(i)\ \beta_8=0,3|a_1a_7, \ 1)\ \textrm { and } \ 5)\ \textrm {  in  Table } \ 13$$
$$ (ii)\ \beta_8\in \{1,5\}, 3|a_2a_5a_8,5|a_3a_8, 10) \ \textrm { and }\ 12)\ \textrm { in Table } 13.$$
 
Let (ii) hold. Consider 10). We have  
$a_6=1, a_2=3$ and \eqref{chi} is not valid with $j_{7^2}=1$ and $(j_1, j_2)=(2, 6)$ and hence 
excluded. The possibility $12)$ implies that $a_j=8-j$ for $1\leq j<8, j\neq 4$. The equation $[2,5,6]$ gives rise to an equation of the form
$x^7_6+x^7_2=2x^7_5$ which is a contradiction by Lemma \ref{eq1} with $p=2$.

Let $(i)$ hold.  Possibility $5)$ gives 
$$a_5=1, a_2=2, a_7\in \{3, 3^6\},\ell=7,q=7^2,j_q=8.$$
By \eqref{chi} with $(j_1,j_2)=(2,5)$ we get $\lambda_q=1.$ Then we check that
\eqref{chi} is not satisfied with $(j_1,j_2)=(6,7).$ Thus $5)$ is excluded.  
Next we take Possibility $1).$ Then 
$$a_2=1, a_5=4$$ and it is easy to check that 
$$a_j=j-1\ \textrm { for }\ 2\leq j\leq 7, j\neq 4\ \textrm { and }\ a_1=2^3\cdot3^6\cdot5^6\cdot7^6.$$ 
We use Lemma \ref{l^2} with $\ell=7, q=29$. 
Suppose $29|d.$ Then taking $(j_1,j_2)=(2,3),$ we see that \eqref{q|d} is not satisfied.
Suppose $29\nmid d.$ Let $j_{29}\not\in \{2,3\}.$ Then \eqref{chi} with $(j_1,j_2)=(2,3)$ implies that
$$j_{29}\in \{ 1,7,9,12\}.$$
Suppose $j_{29}\in \{1,2,3,7,9,12\}.$ Then taking $(j_1,j_2)=(5,7),$ we get $j_{29}=1.$ Thus $29|x_1.$  
Again, using Lemma \ref{l^2} with $\ell=7, q=43$ as above,we obtain
$43|x_1$. Thus $29\cdot 43|x_1$. Now $[1,2,6]$ gives rise to the equation
$$x_2^7 -x_6^7=2^5\cdot 3^6\cdot5^5\cdot7^6(x_1)^7.$$ 
From Lemma \ref{n-n-n-lemma},  we have a contradiction since $29\cdot 43|x_1$. 

 \section{\bf Proof of Theorems \ref{Th2} and  \ref{Th1-3-8} }\label{PfTh1}
 
In this section, we prove Theorems \ref{Th2} and \ref{Th1-3-8}. Let $4 \leq k\leq 8.$ From the sections $4-7$, it suffices to 
consider $\ell=2$. 
Throughout this section,
we consider the following equation for rationals $x, y$ and $1< i < k.$
\begin{equation}\label{case-l=2}
(x+1)\cdots(x+i-1)(x+i+1)\cdots(x+k)=y^2, \quad 4 \leq k\leq 8, 
\end{equation}
As we mentioned before,by symmetry, it is enough to consider $i\leq \frac{k+1}{2}$.

We state the following lemma which will be essential for the proof.
\begin{lemma}\label{HyperEll}
Let $b$ be an integer with $P(b)\leq 5$. We consider the hyperelliptic curve
\begin{equation}\label{HypEll}
b(X+1)\cdots\widehat{(X+j)}\cdots(X+6)=Y^2, \quad 1\leq j \leq 3
\end{equation}
where $(X+j)$ is missing term in the product on the left hand side. Then we have the following:
\begin{enumerate}[(i)]
\item The hyperelliptic curves \eqref{HypEll} have genus $2$.
\item The trivial rational points(i.e., with $Y=0$) on \eqref{HypEll} appears for 
\begin{align*}
X= \begin{cases}
-1,\cdots,-(j-1),-(j+1),-6 & {\rm if} \ j\geq2;\\
-2,\cdots,-6 & {\rm if} \ j=1.
\end{cases}
\end{align*}
\item The non-trivial rational points of \eqref{HypEll} are explicitly determined and they are given by the tuple $(X,j,b)$ in the set\\
$\Omega=\{(0,1,5),(-8,1,-5),(-1,1,30),(-7,1,-30),(-7,2,-1),\\(-9,2,-5),(-2,2,-6),(0,2,10),(-3,3,3),(-7,3,-5),(4,3,6),\\(0,3,15),(-10,3,-15)\}$
\end{enumerate}
\end{lemma}

\begin{proof}
From equation\eqref{HypEll}, we claim that it is enough to find the rational points in the following three hyperelliptic curves
\begin{equation}\label{hyp1}
b(X-2)(X-1)X(X+1)(X+2)=Y^2
\end{equation}
\begin{equation} \label{hyp2}
b(X-3)(X-1)X(X+1)(X+2)=Y^2
\end{equation}
\begin{equation}\label{hyp3}
b(X-3)(X-2)X(X+1)(X+2)=Y^2
\end{equation}
where $P(b)\leq 5.$\\

\noindent
The genus of the curves \eqref{hyp1},\eqref{hyp2}, and \eqref{hyp3} are $2$ and the we can find the corresponding Jacobians($J$). We compute the rank bounds($r$) for the corresponding Jacobians($J$) using a quick program in Magma, we find that $r\leq 1.$ We also know the trivial rational points on these curves, hence, in principle, we can use the method of Chabauty.\\

\noindent
The trivial case, when the Jacobian($J$) has Mordell-Weil rank $r=0$, we find that there are no non-trivial rational points in Magma using the implemented code {\it{Chabauty0(J)}}.\\

\noindent
The non-trivial case, i.e., the Jacobian($J$) has Mordell-Weil rank $r=1$. We use the Chabauty method for genus $2$ curves implemented by Bruin and Stoll \cite{BrSt} in Magma that computes all the rational points. We use primes $7,11,23$ for the Chabauty arguments and we precisely find all the rational points of \eqref{hyp1},\eqref{hyp2}, and \eqref{hyp3} and they are given in the Lemma \ref{HyperEll} (ii) and (iii). This completes the proof.
\end{proof}

Let $4 \leq k\leq 8$ and $\ell=2.$

\vskip 2mm
\noindent
{\bf Case 1. $k=4.$}  Then $i=2$ and we get
$$(X+1)(X+3)(X+4)=Y^2.$$
This is of rank 0 and the torsion points belong to $\mathbb Z/2\mathbb Z\times
\mathbb Z/2\mathbb Z$ leading to only solutions with $y=0.$
\vskip 2mm
\noindent
{\bf Case 2. $k=5.$} Then $i=2,3.$ For $i=2,$ we have the equation
\begin{equation}\label{bad_boy}
(n+s)(n+3s)(n+4s)(n+5s)=y^2
\end{equation}
Taking $X= \frac{n}{s}, Y=\frac{y}{s^2}$ leads to the curve $X(X+2)(X+3)(X+4)=Y^2$ which is birationally equivalent to the elliptic curve $$F: y^2= x^3 + 8x^2 + 12x.$$ This elliptic curve has rank $1$ and the Mordell Weil group of $F$ is isomorphic to $\Z\times \Z_2 \times \Z_2.$ The non-trivial rational points of $F$, leads to infinitely many solutions to the equation \eqref{bad_boy}.
\vskip 2mm
\noindent
For $i=3,$ this is an elliptic curve of rank $0$, hence the only rational points appears for $X=-1,-2,-4$ which corresponds to the torsion points of the curve.\\
\vskip 2mm
\noindent
{\bf Case 3. $k=6.$} Then $i=2,3.$ It clearly follows from Lemma\ref{HyperEll} that there are no non-trivial rational points.
\vskip 5mm
\noindent
{\bf Case 4.} Let $k=7.$ Then $i\in\{2,3,4\}.$ 
\vskip 5mm
Let $i=2.$
Then 
\begin{equation}\label{k=7,l=2,i=2}
(x+1)(x+3)(x+4)(x+5)(x+6)(x+7)=y^2
\end{equation}
In this case, we use the method of Chabauty for genus $2$ curves due to Bruin and Stoll \cite{BrSt}. The hyperelliptic curve \eqref{k=7,l=2,i=2} is of genus $2$ and the Mordell Weil group of the Jacobian, J of the curve \eqref{k=7,l=2,i=2} is given by ${\mathbb{Z}_2}^4\times\mathbb{Z}$ and the generator of the free part is given by the $(x + 1, -x^3 - 1, 2)$. We can use the Chabauty method for genus $2$ curves implemented by Bruin and Stoll \cite{BrSt} in Magma that computes all the rational points. A quick search in Magma finds that the rational points occurs at $x=1,3,4,5,6,7; y=0$ and $x=\frac{-37}{7}; y= \pm \frac{720}{7^3}$. Hence the only non-trivial rational points of \eqref{k=7,l=2,i=2} is $(x,y)=\Big(\frac{-37}{7}, \pm \frac{720}{7^3}\Big)$.

\vskip 5mm
\noindent
Let  $i=3, 4.$ Then we have 

\begin{equation}\label{case-l=2-k=7-i=3}
(x+1)(x+2)(x+4)(x+5)(x+6)(x+7)=y^2
\end{equation}
\begin{equation}\label{case-l=2-k=7-i=4}
(x+1)(x+2)(x+4)(x+5)(x+6)(x+7)=y^2
\end{equation}
\noindent
We implement the same techniques as we did for the case $(k,\ell,i)=(7,2,2)$. We find that there are no non-trivial rational points in this cases.
\vskip 5mm
\noindent
{\bf Case 5.} Let $k=8.$ Then $i\in \{2,3,4\}.$ We have the equation
\begin{equation}\label{eqn-k=8,l=}
(n+s)\cdots\widehat{(n+is)}(n+8s)=y^2,
\end{equation}
where $(n+is)$ is the missing from product in the left hand side.\\
\noindent
Note that here $s=d^2$ for some $d.$ 
\noindent
We can assume that $7\nmid a_2, a_3,\cdots,a_7$, where we fix a convention that if $(n+is)$ is the missing term, then the corresponding $a_i=1.$ This boils down to finding integer points for the equation
\begin{equation}\label{HypEll-Int}
(x+2s)\cdots\widehat{(x+js)}\cdots(x+7s)=by^2, \quad 2\leq j \leq 4, P(b)\leq 5.
\end{equation}
Since $s=d^2$ and $s\neq 0$, we can claim that it is equivalent to finding rational points to the hyperelliptic curve \eqref{HypEll}. By Lemma \ref{HypEll} we check that none of those non-trivial rational points gives rise to a a non-trivial rational points in our case. This completes the proof.
\qed
\vskip 5mm
\section*{Acknowledgement}
 Authors like to thank the referee for careful reading and valuable suggestions. It leads to an improvement of the article.
Pranabesh wants to thank Prof. Cameron Stewart and the Department of Pure Mathematics at the University of Waterloo for the very nice working  environment where part of this project was accomplished. Pranabesh also acknowledges the continuous support of Prof. Karl Dilcher and the postdoctoral grant from the Atlantic Association for Research in the Mathematical Sciences (AARMS), Canada.
Saradha likes to thank Indian National Science Academy for awarding Senior Scientist fellowship under which this work was done. She also thanks DAE-Center for Excellence in Basic Sciences, Mumbai University for providing facilities to carry out this work. Divyum acknowledges the support of the DST-SERB SRG Grant (SRG/2021/000773) and the OPERA award of BITS Pilani.


\end{document}